\documentclass[10pt]{article} 
\usepackage{latexsym}  
\usepackage[all]{xy} 
\usepackage{amsmath} 
\usepackage{amsfonts}

\def\abstractname{R\'esum\'e} 
\newtheorem{thm}{Th\'eor\`eme}[section] 
\newtheorem{prop}[thm]{Proposition} 
\newtheorem{lem}[thm]{Lemme} 
\newtheorem{df}[thm]{D\'efinition} 
\newtheorem{cor}[thm]{Corollaire} 
 
\begin{document}
  
\title{Vers une interpr\'etation galoisienne de la th\'eorie de l'homotopie} 
\author{B. To\"en} 
 
\maketitle 
 
\def\abstractname{Abstract} 
 
\begin{abstract} 
Given any $CW$ complex $X$, and $x \in X$, it is well known that   
$\pi_{1}(X,x)\simeq Aut(\omega_{x}^{0})$, where $\omega_{x}^{0}$ is the functor which associates to each locally  
constant sheaf on $X$ its fibre at $x$. The purpose of the present work 
is to generalize this formula to higher homotopy. For this we introduce the $1$-Segal  
category of locally constant \ ($\infty$-)stacks  
on $X$, and we prove that the $H_{\infty}$-space of endomorphisms of its fibre functor  
at $x$ is equivalent to the loop space $\Omega_{x}X$. 
\end{abstract} 
 
\def\abstractname{R\'esum\'e} 
 
\begin{abstract} 
\'Etant donn\'e un $CW$ complexe $X$, et $x \in X$, il est bien connu que 
$\pi_{1}(X,x)\simeq Aut(\omega_{x}^{0})$, o\`u $\omega_{x}^{0}$ est le foncteur qui associe  
\`a tout faisceau localement constant sur $X$ sa fibre 
en $x$. Le but du pr\'esent travail est de g\'en\'eraliser cet isomorphisme \`a  
l'homotopie sup\'erieure. Pour cela nous introduisons 
la $1$-cat\'egorie de Segal des ($\infty$-)champs localement constants sur $X$,  
et nous montrons que le $H_{\infty}$-espace  
des endomorphismes de son foncteur fibre en $x$ est \'equivalent  
\`a l'espace des lacets $\Omega_{x}X$.   
\end{abstract} 
 
\newpage 
 
\tableofcontents 
 
\newpage 
 
\setcounter{section}{-1} 
 
\begin{section}{Introduction} 
 
Pour tout $CW$ complexe point\'e et connexe $(X,x)$, il est bien connu que la monodromie 
induit un isomorphisme naturel d\'ecouvert par A. Grothendieck 
$$\pi_{1}(X,x)\simeq End(\omega^{0}_{x}) \quad (\text{isomorphisme de mono\"{\i}des}),$$ 
o\`u $\omega^{0}_{x}$ est le foncteur qui associe \`a tout faisceau localement  
constant sur $X$ sa fibre en $x$ (en particulier tout 
endomorphisme de $\omega^{0}_{x}$ est en r\'ealit\'e un automorphisme). Ainsi, le groupe fondamental  
de l'espace topologique point\'e $(X,x)$ peut \^etre reconstruit 
\`a l'aide de la cat\'egorie des faisceaux localement constants 
sur $X$ munie de son foncteur fibre. Le but de ce travail est  
de g\'en\'eraliser ce r\'esultat \`a l'homotopie sup\'erieure.  
Pour cela nous utiliserons des analogues sup\'erieurs des notions de groupe fondamental,  
de faisceau localement constant, de cat\'egorie et de mono\"{\i}de des endomorphismes 
d'un foncteur. Ces notions sup\'erieures seront  
respectivement celles de $H_{\infty}$-espace des lacets, de \textit{($\infty$-)champs localement constants} 
(voir Def. \ref{d7}), 
de \textit{$1$-cat\'egories de Segal} (voir Def. \ref{d4}), et de  
\textit{$H_{\infty}$-espace des  
endomorphismes d\'eriv\'ees d'un foncteur simplicial} 
(voir Def. \ref{d6}), dont nous donnerons un apper\c{c}u dans la suite
de cette introduction. Notre r\'esultat principal (voir Cor. \ref{c5}) est alors la construction  
d'une \'equivalence faible naturelle de $H_{\infty}$-espaces (voir Def. \ref{d3})
$$\Omega_{x}X \simeq \mathbb{R}End(\omega_{x}).$$ 
Ainsi, comme $X$ est faiblement \'equivalent  
au classifiant $B\Omega_{x}X$, le type d'homotopie d'un  
$CW$ complexe point\'e et connexe $(X,x)$ peut \^etre reconstruit  
\`a l'aide de la $1$-cat\'egorie de Segal  
des $(\infty-)$champs localement constants sur $X$ munie de son foncteur fibre.  
L'\'equivalence faible $\Omega_{x}X \simeq \mathbb{R}End(\omega_{x})$ m\'erite 
sans doutes le nom d'\textit{$\infty$-monodromie universelle}  
(voir la remarque \`a la fin du $\S 3$). \\ 
 
\begin{center} \textit{L'isomorphisme  
$\pi_{1}(X,x)\simeq End(\omega_{x}^{0})$} \end{center}  
 
Commen\c{c}ons par rappeler bri\`evement une construction de l'isomorphisme 
$\pi_{1}(X,x)\simeq End(\omega_{x}^{0})$ cit\'e plus haut, qui nous servira de mod\`ele 
tout au long de ce travail.  
 
Consid\'erons $Loc_{0}(X)$ la 
cat\'egorie des faisceaux d'ensembles localement constants sur $X$. On dispose  
du foncteur fibre en $x$, $\omega^{0}_{x} : Loc_{0}(X) \longrightarrow Ens$. 
La construction qui associe \`a un faisceau localement constant sa monodromie induit  
une \'equivalence entre la cat\'egorie $Loc_{0}(X)$, et  
la cat\'egorie $\pi_{1}(X,x)-Ens$, des ensembles munis 
d'une action (disons \`a gauche)  
du groupe $\pi_{1}(X,x)$. A travers cette \'equivalence le foncteur fibre $\omega^{0}_{x}$ est  
transform\'e en le foncteur d'oubli  
$\omega^{0} : \pi_{1}(X,x)-Ens \longrightarrow Ens$, et on a 
donc $End(\omega^{0}_{x})\simeq End(\omega^{0})$.  
Enfin, le foncteur $\omega^{0}$ est co-repr\'esent\'e par 
l'objet $E$ de $\pi_{1}(X,x)-Ens$ d\'efini par l'action de $\pi_{1}(X,x)$ sur lui-m\^eme par translations.  
Le lemme de Yoneda implique alors l'existence d'isomorphismes naturels 
$$End(\omega^{0}_{x})\simeq End(\omega^{0})\simeq End(E).$$  
Enfin, il est imm\'ediat de constater que l'application qui envoie 
un endomorphisme $\phi$ de $E$ sur $\phi(e)^{-1}\in \pi_{1}(X,x)$ 
induit un isomorphisme de mono\"{\i}des $End(E) \simeq \pi_{1}(X,x)$. En conclusion,  
il existe une chaine d'isomorphismes naturels de mono\"{\i}des 
$$End(\omega^{0}_{x})\simeq End(\omega^{0})\simeq End(E)\simeq \pi_{1}(X,x).$$ 
Si l'on explicite ces isomorphismes, on peut comprendre l'isomorphisme 
$\pi_{1}(X,x)\simeq End(\omega^{0}_{x})$ de la fa\c{c}on suivante. 
Pour un lacet $\gamma \in \pi_{1}(X,x)$, et un faisceau localement constant 
$F \in Loc_{0}(X)$, la monodromie de $F$ le long de $\gamma$ induit un 
automorphisme $\gamma_{F} : \omega^{0}_{x}(F) \simeq \omega^{0}_{x}(F)$. Il est facile de voir 
que cet automorphisme est fonctoriel en $F$, et d\'efinit donc 
un automorphisme $[\gamma]$ du foncteur $\omega^{0}_{x}$. On peut 
v\'erifier que l'isomorphisme $\pi_{1}(X,x)\simeq End(\omega^{0}_{x})$ d\'ecrit ci-dessus 
est induit par la correspondance $\gamma \mapsto [\gamma]$.  
Notre approche est de proc\'eder par analogie et de 
g\'en\'eraliser pas \`a pas la d\'emarche pr\'ec\'edente \`a l'homotopie sup\'erieure. 
 
\begin{center} \textit{Champs localement constants} \end{center} 
 
Pour expliquer notre point de vue sur la notion de champs rappelons une construction 
(non conventionnelle) du topos de l'espace $X$  
(i.e. d'une cat\'egorie qui est naturellement \'equivalente  
\`a la cat\'egorie des faisceaux sur $X$) qui n'utilise pas directement la notion de faisceaux. 
Pour cela, soit $Pr(X)$ la cat\'egorie des pr\'efaisceaux d'ensembles sur l'espace topologique 
$X$. Dans cette cat\'egorie on consid\`ere l'ensemble $W$ des morphismes qui induisent 
des isomorphismes fibre \`a fibre, et on forme la cat\'egorie $W^{-1}Pr(X)$, obtenue 
\`a partir de $Pr(X)$ en inversant formellement les morphismes de $W$. On peut alors 
v\'erifier que $W^{-1}Pr(X)$ est naturellement \'equivalente \`a la cat\'egorie des faisceaux sur $X$.  
Il faut remarquer que les objets de $W^{-1}Pr(X)$ sont les pr\'efaisceaux sur 
$X$, mais ses ensembles de morphismes sont en r\'ealit\'e isomorphes aux ensembles 
de morphismes entre faisceaux associ\'es. Aussi surprenant que cela 
puisse para\^{\i}tre, cette construction montre 
qu'il n'est pas n\'ecessaire de conna\^{\i}tre la notion de faisceaux pour pouvoir 
parler de la cat\'egorie des faisceaux sur $X$.  
 
Nous d\'efinirons la cat\'egorie simpliciale\footnote{Tout au long  
de ce travail, l'expression  
\textit{cat\'egorie simpliciale} fera 
r\'ef\'erence \`a la notion de cat\'egorie enrichie dans les  
ensembles simpliciaux. Il s'agit donc 
de la notion utilis\'ee dans \cite{dk}, et non de la notion  
plus g\'en\'erale d'objet simplicial dans 
la cat\'egorie des cat\'egories.} des champs sur $X$ en s'appuyant sur 
le m\^eme principe, \`a savoir qu'il n'est pas n\'ecessaire de conna\^{\i}tre la notion 
de \textit{champs} pour pouvoir parler de la \textit{cat\'egorie simpliciale des champs sur $X$}.  
Ainsi, il ne s'agit pas de donner une d\'efinition pr\'ecise de \textit{champs}, mais 
de construire une \textit{cat\'egorie simpliciale des champs sur $X$} qui est raisonnable. 
Les objets de cette cat\'egorie n'auront pas grande importance par la suite, et c'est  
la cat\'egorie dans son ensemble \`a laquelle nous nous int\'eresserons. 
 
Soit $SPr(X)$ la cat\'egorie des pr\'efaisceaux d'ensembles simpliciaux sur l'espace topologique $X$.  
Dans cette cat\'egorie on consid\`ere l'ensemble $W$ des morphismes qui induisent 
des \'equivalences faibles fibre \`a fibre. La cat\'egorie simpliciale des champs sur $X$  
est par d\'efinition la localis\'ee simpliciale de Dwyer-Kan $LSPr(X):=L(SPr(X),W)$, obtenue \`a partir 
de $SPr(X)$ en inversant formellement les morphismes de $W$ (voir \cite{dk} et Def. \ref{d6''}).  
On d\'efinit ensuite la cat\'egorie simpliciale des champs localement 
constants sur $X$ comme la sous-cat\'egorie simpliciale pleine de $LSPr(X)$ form\'ee des  
objets qui sont localement \'equivalents \`a des champs constants (voir Def. \ref{d7}). 
C'est cette derni\`ere cat\'egorie simpliciale qui nous int\'eressera par la suite, et qui  
sera not\'ee $Loc(X)$. Par d\'efinition $Loc(*)$ sera simplement not\'ee $Top$, et est la  
localis\'ee simpliciale de Dwyer-Kan de la cat\'egorie des ensembles simpliciaux  
le long des \'equivalences faibles. Enfin, il est clair que la construction $X \mapsto Loc(X)$ est  
fonctorielle en $X$, et donc le morphisme $* \longrightarrow X$ correspondant au point de 
base $x \in X$ induit un foncteur simplicial 
$$\omega_{x} : Loc(X) \longrightarrow Loc(*)=:Top,$$ 
qui sera notre foncteur fibre.  \\ 
 
\begin{center} \textit{Cat\'egories simpliciales, $1$-cat\'egories de Segal et  
$\mathbb{R}End(\omega_{x})$} \end{center} 
 
Nous venons de d\'efinir la cat\'egorie simpliciale $Loc(X)$ des champs localement constants sur $X$, et 
son foncteur fibre $\omega_{x} : Loc(X) \longrightarrow Top$. Notre probl\`eme est maintenant de 
d\'efinir de fa\c{c}on convenable un ensemble simplicial des endomorphismes de $\omega_{x}$. 
Bien entendu, il existe un ensemble simplicial na\"{\i}f $End(\omega_{x})$, des endomorphismes  
simpliciaux de $\omega_{x}$. Cependant, $End(\omega_{x})$ ne poss\`ede pas \textit{le bon type 
d'homotopie}, et il est fort probable que notre r\'esultat principal soit faux  
si on utilise cette d\'efinition. Une fa\c{c}on de r\'esoudre ce probl\`eme est  
de trouver une version \textit{d\'eriv\'ee} de la construction $\omega_{x} \mapsto End(\omega_{x})$. 
Pour cela, nous avons choisi d'utiliser la 
th\'eorie des $1$-cat\'egories de Segal et de ses $Hom$ internes (voir \cite{s1,s3}).  
La notion de $1$-cat\'egorie de Segal 
est une notion affaiblie de celle de cat\'egorie simpliciale, exactement de la m\^eme mani\`ere que 
la notion de $H_{\infty}$-espace est une notion affaiblie de celle de groupe simplicial.  
Les cat\'egories simpliciales s'identifient de plus aux $1$-cat\'egories de Segal 
\textit{strictes}, et il est vrai que toute $1$-cat\'egorie de Segal est \'equivalente  
\`a une cat\'egorie simpliciale (ceci est une g\'en\'eralisation du fait que 
tout $H_{\infty}$-espace est \'equivalent \`a un groupe simplicial). En r\'ealit\'e, 
les th\'eories homotopiques des cat\'egories simpliciales et des cat\'egories de Segal 
sont \'equivalentes (voir le paragraphe \textit{Strictification} de \cite{s3}). 
Le lecteur pourra donc en premi\`ere approximation remplacer 
mentalement les $1$-cat\'egories de Segal par des cat\'egories simpliciales. 
Cependant, un avantage de nature technique des $1$-cat\'egories de Segal est que 
l'on peut d\'efinir des $Hom$ internes en un sens d\'eriv\'e (voir $\S 1.2$),  
qui seront not\'es  
par la suite $\mathbb{R}\underline{Hom}$\footnote{Lorsque $A$ et $B$ sont deux cat\'egories 
simpliciales, les objets de $\mathbb{R}\underline{Hom}(A,B)$ peuvent \^etre 
l\'egitimement per\c{c}us comme des foncteurs simpliciaux \textit{faibles} entre $A$ et $B$}. 
 
Revenons \`a notre foncteur fibre $\omega_{x} : Loc(X) \longrightarrow Top$, qui  
sera alors consid\'er\'e comme un objet dans la $1$-cat\'egorie de Segal 
$\mathbb{R}\underline{Hom}(Loc(X),Top)$. L'ensemble simplicial 
des endomorphismes de cet objet sera not\'e $\mathbb{R}End(\omega_{x})$. L'ensemble simplicial 
$\mathbb{R}End(\omega_{x})$ vient avec une loi de composition induite  
par la composition des endomorphismes qui le 
munit d'une structure naturelle de $H_{\infty}$-espace. \\ 
 
\begin{center} \textit{La preuve dans ses grandes lignes} \end{center} 
 
Nous sommes maintenant en mesure de suivre pas \`a pas la preuve de l'isomorphisme 
$\pi_{1}(X,x)\simeq End(\omega^{0}_{x})$ pr\'esent\'ee pr\'ec\'edemment. 
 
Nous commencerons par g\'en\'eraliser la classification des faisceaux localement constants par leur 
monodromie au cas des champs localement constants d\'efinis plus haut.  
Cette classification se fera en deux grandes \'etapes. 
Tout d'abord,  
nous montrerons que la $1$-cat\'egorie de Segal $Loc(X)$, des champs localement 
constants sur $X$ est \'equivalente \`a la $1$-cat\'egorie de Segal $Fib(S(X))$, des fibrations sur  
l'ensemble simplicial des simplexes singuliers de $X$ (voir Thm. \ref{p3}). Cette \'equivalence peut \^etre 
comprise comme un analogue du fait que la cat\'egorie des faisceaux localement constants est \'equivalente 
\`a celles des rev\^etements. 
Par la suite, nous  
noterons $G$ le groupe simplicial des lacets sur $S(X)$ de base $x$ (i.e. un groupe 
simplicial $G$ faiblement \'equivalent au $H_{\infty}$-espace $\Omega_{x}X$), et nous consid\'ererons 
la $1$-cat\'egorie de Segal $G-Top$ des ensembles simpliciaux munis d'une action 
de $G$ (voir Def. \ref{d7'}). 
Nous montrerons alors que la $1$-cat\'egorie de Segal $Fib(S(X))$ est \'equivalente 
\`a $G-Top$, la conclusion \'etant que $Loc(X)$ est elle m\^eme \'equivalente \`a $G-Top$ (voir Thm. \ref{p5}).  
Cette derni\`ere \'equivalence est l'analogue final de la classification des faisceaux  
localement constants par leur monodromie (i.e. de l'\'equivalence entre faisceaux localement constants 
et $\pi_{1}(X,x)$-ensembles). De plus, le foncteur fibre $\omega_{x} : Loc(X) \longrightarrow Top$ sera alors 
transform\'e en le foncteur d'oubli $\omega_{G} : G-Top \longrightarrow Top$, ce qui 
implique que le $H_{\infty}$-espace $\mathbb{R}End(\omega_{x})$ est faiblement \'equivalent 
\`a $\mathbb{R}End(\omega_{G})$. 
 
Pour achever la d\'emonstration, nous remarquerons que le foncteur $\omega_{G}$ est co-repr\'esentable 
(au sens des $1$-cat\'egories de Segal) par l'objet de $G-Top$ d\'efini par l'action  
de $G$ sur lui-m\^eme. La version $1$-cat\'egorie de Segal du lemme de Yoneda impliquera alors 
l'existence d'une \'equivalence faible d'ensembles simpliciaux 
$$\mathbb{R}End(\omega_{x})\simeq \mathbb{R}End(\omega_{G}) \simeq G\simeq \Omega_{x}X,$$ 
ce qui est bien le r\'esultat annonc\'e.  
 
Remarquons pour terminer l'\'equivalence faible de notre th\'eor\`eme 
$\Omega_{x}X \simeq \mathbb{R}End(\omega_{x})$, 
peut alors \^etre comprise comme induite par l'action de la monodromie des lacets de $X$ 
sur les fibres des champs localement constants. Ainsi, intuitivement parlant,  
un lacet $\gamma$ sur $X$ correspond dans $\mathbb{R}End(\omega_{x})$  
\`a l'endomorphisme de $\omega_{x}$ donn\'e par la 
monodromie $[\gamma]$ le long de $\gamma$, une homotopie entre deux lacets  
$\gamma$ et $\gamma'$ correspond \`a l'homotopie entre les endomorphismes 
$[\gamma]$ et $[\gamma']$ donn\'ee par une certaine \textit{monodromie  
sur un pav\'e de dimension $2$} \dots et ainsi de suite, tout cela de fa\c{c}on compatible 
avec la composition des lacets et la composition des endomorphismes. \\ 
 
\begin{center} \textit{Organisation du travail} \end{center}  
 
Dans une premi\`ere partie nous avons rappel\'e un certain nombre de d\'efinitions et  
de r\'esultats pr\'eliminaires qui nous 
seront utiles tout au long de l'article.  
Il comporte en particulier les notations et d\'efinitions concernant 
les $H_{\infty}$-espaces, les $1$-cat\'egories de Segal et la localisation  
de Dwyer-Kan.  
 
La seconde partie est consacr\'ee \`a la d\'emonstration des r\'esultats de classification 
des champs localement constants, et est le coeur du travail. Les preuves utilisent 
le formalisme des cat\'egories de mod\`eles et des adjonctions de Quillen afin de construire  
des \'equivalences entre les diff\'erentes $1$-cat\'egories de Segal qui entrent en jeu. 
 
Enfin, dans une derni\`ere partie nous rassemblons les r\'esultats pr\'ec\'edemment d\'emontr\'es et  
en d\'eduisons le th\'eor\`eme principal. \\ 
 
\begin{center} \textit{Relations avec d'autres travaux} \end{center} 
 
Dans la lettre \`a L. Breen, dat\'ee du $17/02/1975$ et 
que l'on trouve dans \cite{gr}, A. Grothendieck conjecture 
plusieurs relations entre la th\'eorie des $n$-champs en groupo\"{\i}des et la th\'eorie de 
l'homotopie. Bien qu'utilisant un language tr\`es diff\'erent,  
ces consid\'erations ont \'et\'e une source d'inspiration 
lors de l'\'elaboration de ce travail (voir la remarque qui suit le corollaire \ref{c5}).  
 
Signalons aussi que le corollaire \ref{c5} avait \'et\'e \'enonc\'e sans 
d\'emonstration dans l'introduction de \cite{to}, et sert d\'ej\`a depuis quelques temps 
de point de d\'epart au formalisme Tannakien sup\'erieur qui y est propos\'e.  
 
\newpage 
 
\begin{center} \textit{Notations et conventions} \end{center} 
 
Pour tout univers $\mathbb{U}$ nous noterons $Ens_{\mathbb{U}}$ (resp. $SEns_{\mathbb{U}}$, resp.  
$Es_{\mathbb{U}}$,  resp. \dots) la cat\'egorie des ensembles (resp. des ensembles 
simpliciaux, resp. des espaces topologiques, resp. \dots) appartenant \`a $\mathbb{U}$.  
Lorsque le contexte ne n\'ecessitera qu'un unique univers, nous 
omettrons de le mentionner, et nous \'ecrirons simplement $Ens$, $SEns$, $Es$, \dots. 
 
Nous utiliserons l'abr\'eviation \textit{cmf} pour signifier \textit{cat\'egorie  
de mod\`eles ferm\'ee}, et pour nous cette expression fera toujours r\'ef\'erence 
\`a la notion expos\'ee dans \cite[$\S 1.1.3$]{ho} (en particulier, nos cmf poss\`ederons toujours 
tout type de limites et colimites ainsi que des factorisations fonctorielles). 
L'expression \textit{\'equivalence} sera synonyme d'\textit{\'equivalence faible}, et fera 
r\'ef\'erence \`a une structure de cat\'egorie de mod\`eles.  
 
Pour une cmf $M$,  et $X$ un de ses objets, nous noterons $M/X$ (resp. $X/M$)  
la cmf des objets au-dessus (resp. au-dessous) de $X$. Nous noterons 
aussi $*$ l'objet final d'une cmf $M$ (choisit une bonne fois pour toute),  
et $*/M$ la cmf des objets point\'es de $M$. Nous utiliserons g\'en\'eralement  
la notation $Hom_{*}$ pour d\'esigner les morphismes d'objets point\'es. 
La cat\'egorie homotopique d'une cmf $M$ sera not\'ee $Ho(M)$.  
 
Un foncteur $F : M \longrightarrow N$ entre deux cmf sera de Quillen  
\`a droite (resp. \`a gauche) si c'est un adjoint \`a droite (resp. \`a gauche) 
et s'il pr\'es\`erve les fibrations et les fibrations triviales (resp. 
les cofibrations et les cofibrations triviales). 
Si $F : M \longrightarrow N$ est un foncteur de Quillen \`a droite entre 
deux cmf, d'adjoint \`a gauche $G : N \longrightarrow M$, nous noterons  
$$\mathbb{R}F : Ho(M) \longrightarrow Ho(N) \qquad \mathbb{L}G : Ho(N) \longrightarrow Ho(M)$$ 
leurs foncteurs d\'eriv\'es \`a droite et \`a gauche (\cite[$\S 1.3$]{ho}).  
 
Nous utiliserons aussi la notion de cmf simpliciale de \cite[$\S 4.2.18$]{ho}. 
Si $N$ est une cmf simpliciale, et si $\underline{Hom}$ d\'esigne ses ensembles  
simpliciaux de morphismes, nous noterons  
$\mathbb{R}\underline{Hom}$ ses $Hom$ d\'eriv\'es \`a valeurs dans $Ho(SEns)$ (voir \cite[$\S 4.3.4$]{ho}). 
Rappelons qu'ils sont d\'efinis pour deux objets $X$ et $Y$ de $Ho(M)$, par 
$\mathbb{R}\underline{Hom}(X,Y):=\underline{Hom}(QX,RY) \in Ho(SEns)$,  
o\`u $QX$ (resp. $RY$) est un mod\`ele cofibrant (resp. fibrant) 
de $X$ (resp. de $Y$).  
Ceux-ci munissent la cat\'egorie $Ho(M)$ d'une structure 
de $Ho(SEns)$-modules ferm\'e au sens de \cite[$\S 4.1.13$]{ho} (i.e. 
$Ho(M)$ est enrichie dans $Ho(SEns)$ et poss\`ede des produits externes ainsi que des 
objets exponentiels par des objets de $Ho(SEns)$). 
 
Nous noterons comme il est en l'usage $\Delta$ la cat\'egorie simpliciale  
standard. Ses objets sont les ensembles finis ordonn\'es 
$[m]:=\{0,\dots,m\}$, et ses morphismes les morphismes croissants au sens  
large. Nous supposerons toujours que les univers consid\'er\'es 
contiennent $\Delta$.  
 
Enfin, comme nous l'avons d\'ej\`a signal\'e au cours de cette introduction,  
l'expression \textit{cat\'egorie simpliciale} fera toujours 
r\'ef\'erence \`a la notion de cat\'egorie enrichie dans les ensembles simpliciaux. Il s'agit donc 
de la notion utilis\'ee dans \cite{dk} et \cite[$\S XII$]{dk2} (sous le nom  
de \textit{$S$-cat\'egorie}), et non de  
la notion plus g\'en\'erale d'objet simplicial dans la cat\'egorie des cat\'egories.  
 
\end{section} 
 
\newpage 
 
\begin{section}{Pr\'eliminaires} 
 
Dans ce premier chapitre on rappelle plusieurs d\'efinitions et r\'esultats qui nous 
seront utiles tout au long de ce travail. Il s'agit pour la plupart de notion standard et le lecteur 
est invit\'e \`a consid\'erer cette partie comme une annexe de notations et de r\'ef\'erences. 
 
\begin{subsection}{Espaces topologiques, ensembles simpliciaux, groupes simpliciaux et $\Delta^{o}$-espaces} 
 
Commen\c{c}ons par rappeler que les cat\'egories $Es$ des espaces topologiques et  
$SEns$ des ensembles simpliciaux sont des cmf simpliciales (voir \cite[$\S 4.2.18$]{ho}).  
Rappelons que les \'equiavlences dans $Es$ et $SEns$ sont par d\'efinition les 
\'equivalences faibles (i.e. les morphismes induisant des isomorphismes 
sur les groupes d'homotopie). 
De plus, il existe deux foncteurs (voir \cite[$\S 3.6.7$]{ho}) 
$$|.| : SEns \longrightarrow Es \qquad S : Es \longrightarrow SEns,$$ 
d\'efinissant une \'equivalence de Quillen entre les cmf $Es$ et $SEns$. Notons au  
passage que le foncteur $|.|$ est compatible avec  
la structure simpliciale (i.e. commute avec les produits externes par des ensembles simpliciaux).  
 
\begin{df}\label{d1} 
Un pr\'e-$\Delta^{o}$-espace est un foncteur 
$$\begin{array}{cccc} 
A : & \Delta^{o} & \longrightarrow & SEns\\ 
& [m] & \mapsto & A_{m} 
\end{array}$$ 
tel que $A_{0}=*$.  
 
La cat\'egorie des pr\'e-$\Delta^{o}$-espaces sera not\'ee $Pr\Delta^{o}-SEns$.  
\end{df} 
 
Comme $[0]$ est initial dans $\Delta^{o}$, on peut identifier la cat\'egorie 
$Pr\Delta^{o}-SEns$ avec la cat\'egorie des pr\'efaisceaux 
sur $\Delta-[0]$ \`a valeurs dans $*/SEns$. Comme $*/SEns$ est une cmf  
engendr\'ee par cofibrations on peut appliquer  
\cite[$\S II.6.9$]{gj} et munir $Pr\Delta^{o}-SEns$ d'une structure de cmf.  
Rappelons que les fibrations (resp. les \'equivalences) 
pour cette structure sont les morphismes $f : A \longrightarrow B$ tels que 
pour tout $[m] \in \Delta$, le morphisme induit  
$f_{m} : A_{m} \longrightarrow B_{m}$ soit une fibration (resp. une 
\'equivalence). Rappelons que d'apr\`es nos conventions  
\textit{\'equivalence} signifie \textit{\'equivalence 
faible}. \\ 
 
Soit $[m] \in \Delta$. On dispose de $m$ morphismes  
$$p_{i} : [1] \longrightarrow [m] \qquad 0\leq i< m,$$ 
d\'efinis par $p_{i}(0)=i$ et $p_{i}(1)=i+1$. Ceci permet de d\'efinir pour tout  
$A \in Pr\Delta^{o}-SEns$ et $[m] \in \Delta$ un morphisme 
(dit \textit{morphisme de Segal}), $\prod p_{i}  : A_{m} \longrightarrow A_{1}^{m}$. 
 
\begin{df}\label{d2} 
Un objet $A$ de $Pr\Delta^{o}-SEns$ est appel\'e un $\Delta^{o}$-espace  
si pour tout $[m] \in \Delta$ le morphisme de Segal,  
$A_{m} \longrightarrow A_{1}^{m}$, est une \'equivalence dans $SEns$. 
\end{df} 
 
Pour $X \in SEns$ notons $\pi_{0}(X)$ l'ensemble des ses composantes  
connexes. Comme le foncteur $\pi_{0}(-)$ transforme 
\'equivalences en isomorphismes et commute avec les produits  
directs, pour tout $\Delta^{o}$-espace 
$A$, $\pi_{0}(A_{1})$ est naturellement muni d'une loi de  
composition associative et unitaire. C'est donc un mono\"{\i}de. 
 
\begin{df}\label{d3} 
Un objet $A$ de $Pr\Delta^{o}-SEns$ est appel\'e un $H_{\infty}$-espace si c'est une 
$\Delta^{o}$-espace et si de plus le mono\"{\i}de $\pi_{0}(A_{1})$ est un groupe. 
 
On notera $H_{\infty}-SEns$ la sous-cat\'egorie pleine de $Pr\Delta^{o}-SEns$  
form\'ee des $H_{\infty}$-espaces. 
\end{df} 
 
Remarquer que les notions de $\Delta^{o}$-espaces et de $H_{\infty}$-espaces  
sont invariantes par \'equivalences. Ainsi les propri\'et\'es  
d'\^etre un $\Delta^{o}$-espace ou un $H_{\infty}$-espace  
sont stables par isomorphismes dans la cat\'egorie homotopique $Ho(Pr\Delta^{o}-SEns)$. \\ 
 
Pour un objet $A \in Pr\Delta^{o}-SEns$ on d\'efinit sa diagonale 
$$\begin{array}{cccc} 
d(A) : & \Delta^{o} & \longrightarrow & Ens \\ 
& [m] & \mapsto & (A_{m})([m]), 
\end{array}$$ 
qui est un ensemble simplicial. De plus, comme $A_{0}=*$, $d(A)$ est  
r\'eduit (i.e. $d(A)([0])=*$) et en particulier 
peut-\^etre consid\'er\'e sans ambiguit\'e comme un objet de $*/SEns$ (i.e. un ensemble  
simplicial point\'e). Ceci  
d\'efinit \'evidemment un foncteur 
$$d : Pr\Delta^{o}-SEns \longrightarrow */SEns.$$ 
En composant avec le foncteur de r\'ealisation g\'eom\'etrique  
on obtient un foncteur classifiant 
$$B:=|d(-)| : Pr\Delta^{o}-SEns \longrightarrow */Es.$$ 
Il est bien connu que le foncteur $B$ pr\'eserve les \'equivalences, et  
induit donc un foncteur sur les cat\'egories homotopiques 
$$B : Ho(Pr\Delta^{o}-SEns) \longrightarrow Ho(*/Es).$$ 
 
\begin{thm}{(\cite{se})}\label{t1} 
Le foncteur $B$ est de Quillen \`a gauche d'adjoint \`a droite 
$$\Omega_{*} : */Es \longrightarrow Pr\Delta^{o}-SEns.$$ 
De plus, pour tout $X \in */Es$, $\mathbb{R}\Omega_{*}(X)$ est  
un $H_{\infty}$-espace, et le morphisme d'adjonction 
$$B\mathbb{R}\Omega_{*}X\simeq B\mathbb{R}\Omega_{*}X \longrightarrow X$$ 
est un isomorphisme dans $Ho(*/Es)$ lorsque $X$ est connexe par arcs. 
 
En particulier, le foncteur $B : Ho(H_{\infty}-SEns) \longrightarrow Ho(*/Es)$  
est pleinement fid\`ele et son image essentielle 
est form\'ee des espaces point\'es connexes par arcs. 
\end{thm} 
 
Pour m\'emoire rappelons que $\Omega_{*}(X,x)$ est le $H_{\infty}$-espace  
des lacets de base $x$ sur $X$, et est d\'efini par 
$$\begin{array}{cccc} 
\Omega_{*}(X,x) : & \Delta^{o} & \longrightarrow & SEns \\ 
& [m] & \mapsto & \underline{Hom}_{*}(\Delta^{m}_{*},X) 
\end{array}$$ 
o\`u $\Delta^{m}_{*}$ est obtenu \`a partir de $\Delta^{m}$ en identifiant  
tous ses sommets en un point unique. En d'autres 
termes, $\Omega_{*}(X,x)_{m}$ est l'ensemble simplicial des  
morphismes $\Delta^{m} \longrightarrow X$ qui envoient  
tous les sommets de $\Delta^{m}$ sur le point $x$. \\ 
 
Notons $SGp$ la cat\'egorie des groupes simpliciaux  
(i.e. des objets en groupes dans $SEns$, ou encore des objets simpliciaux 
dans la cat\'egorie des groupes). Rappelons que c'est  
une cmf o\`u les fibrations (resp. les \'equivalences) sont les morphismes induisant 
une fibration (resp. une \'equivalence) sur les  
ensembles simpliciaux sous-jacents (voir \cite[$\S II.5.2$]{gj}). 
 
Il existe une inclusion naturelle 
$$\begin{array}{cccc} 
j : & SGp & \longrightarrow & Pr\Delta^{o}-SEns \\ 
& G & \mapsto & ([m] \mapsto G^{m}), 
\end{array}$$ 
o\`u les faces et les d\'eg\'en\'erescences sont induites par  
les multiplications et les diagonales. 
Le foncteur $j$ identifie $SGp$ avec la sous-cat\'egorie 
pleine de $Pr\Delta^{o}-SEns$ form\'ee des objets $A$ pour lesquels  
les morphismes de Segal sont des isomorphismes,  
et tels que la composition $A_{1}\times A_{1} \longrightarrow A_{1}$  
fasse de $A_{1}$ un groupe simplicial.  
Ainsi, tout groupe simplicial sera consid\'ere implicitement  
comme un objet de $Pr\Delta^{o}-SEns$. De plus, le foncteur 
$j$ poss\`ede un adjoint \`a gauche, $Str : Pr\Delta^{o}-SEns \longrightarrow SGp$,  
qui \`a un pr\'e-$\Delta^{o}$-espace 
$X$ associe le groupe simplicial qu'il engendre.  
Comme $j$ est clairement de Quillen \`a droite, le foncteur 
$Str$ est de Quillen \`a gauche. 
 
\begin{prop}\label{p1} 
Le foncteur $j\simeq \mathbb{R}j : Ho(SGp) \longrightarrow Ho(Pr\Delta^{o}-SEns)$  
est pleinement fid\`ele et son image essentielle 
est form\'ee des $H_{\infty}$-espaces. 
\end{prop} 
 
\textit{Preuve:} En utilisant le th\'eor\`eme \ref{t1}, il suffit 
de montrer que le foncteur $d : Ho(SGp) \longrightarrow Ho(*/SEns)$ 
est pleinement fid\`ele et son image essentielle form\'ee des objets connexes.  
 
Pour $G \in SGp$, notons $EG$ la diagonale de l'ensemble bi-simplicial  
$([m],[p]) \mapsto G_{p}^{m+1}$. C'est donc le nerf du morphisme 
d'ensembles simpliciaux $G \longrightarrow *$, et est donc contractile et fibrant.  
De plus, $EG$ poss\`ede une action naturelle de $G$.  
Plus pr\'ecis\'ement, le groupe $G_{m}$ op\`ere diagonalement sur $EG_{m}=G_{m}^{m+1}$  
par son action \`a gauche. On v\'erifie alors facilement que  
$d(G)$ est naturellement isomorphe au quotient de $EG$ par $G$. On peut alors 
utiliser le th\'eor\`eme \cite[$\S V.3.9$]{gj} pour montrer que  
$d(G)$ est naturellement isomorphe dans $Ho(*/SEns)$ \`a  
$\overline{W}G$ d\'efini dans \cite[$\S V.4$]{gj}. En d'autres  
termes les foncteurs $d : Ho(SGp) \longrightarrow Ho(*/SEns)$ et  
$\overline{W} : Ho(SGp) \longrightarrow Ho(*/SEns)$ sont isomorphes. 
Or, d'apr\`es le corollaire \cite[$\S V.6.4$]{gj},  
le foncteur $\overline{W}$ est pleinement fid\`ele et son image essentielle  
est form\'ee des objets connexes. \hfill $\Box$ \\ 
 
Notons $Ho(*/SEns)_{c}$ la sous-cat\'egorie pleine de  
$Ho(*/SEns)$ form\'ee des objets connexes. On peut construire un inverse  
de l'\'equivalence $d : Ho(SGp) \longrightarrow Ho(*/SEns)_{c}$  
de la fa\c{c}on suivante. On commence par consid\'erer le  
foncteur $d : Pr\Delta^{o}-SEns \longrightarrow */SEns$,  
qui envoie un objet $X \in Pr\Delta^{o}-SEns$ sur sa diagonale.  
Ce foncteur poss\`ede un adjoint \`a droite  
$L_{*} : */SEns \longrightarrow Pr\Delta^{o}-SEns$ d\'efini de la m\^eme fa\c{c}on que 
le foncteur $\Omega_{*}$. En clair, on a pour $(S,s) \in */SEns$ 
$$\begin{array}{cccc} 
L_{*}(S,s) : & \Delta^{o} & \longrightarrow & SEns \\ 
& [m] & \mapsto & \underline{Hom}_{*}(\Delta^{m}_{*},X). 
\end{array}$$ 
Le foncteur $L_{*}$ \'etant clairement de Quillen \`a droite,  
$d$ est un foncteur de Quillen \`a gauche. On obtient donc une adjonction 
sur les cat\'egories homotopiques 
$$\mathbb{L}d : Ho(Pr\Delta^{o}-SEns) \longrightarrow Ho(*/SEns)$$
$$\mathbb{R}L_{*} : Ho(*/SEns) \longrightarrow Ho(Pr\Delta^{o}-SEns).$$ 
On v\'erifie alors ais\'ement \`a l'aide du th\'eor\`eme  
\ref{t1} et de la proposition \ref{p1} que le foncteur  
$$G:=\mathbb{L}Str \circ \mathbb{R}L_{*} : Ho(*/SEns)_{c} \longrightarrow Ho(SGp)$$ 
est un inverse du foncteur $d$.  
 
\end{subsection} 
 
\begin{subsection}{$1$-Cat\'egories de Segal} 
 
La notion de $1$-cat\'egorie (resp. $1$-pr\'ecat\'egorie) de Segal  
est une notion afaiblie de celle de cat\'egorie 
simpliciale, de la m\^eme fa\c{c}on que la notion de  
$\Delta^{o}$-espace (resp. pr\'e-$\Delta^{o}$-espace)  
est une notion afaiblie de celle de mono\"{\i}de simplicial.  
Tout comme pour le cas des cat\'egories 
simpliciales (voir \cite[$\S XII$]{dk2}), il existe une structure de cat\'egorie 
de mod\`eles sur la cat\'egorie des $1$-pr\'ecat\'egorie de  
Segal o\`u les \'equivalences sont  
une certaine g\'en\'eralisation de la notion d'\'equivalence de cat\'egories. De plus,  
on sait que la th\'eorie homotopique des $1$-cat\'egories de Segal est \'equivalente \`a celle 
des cat\'egories simpliciales (voir par exemple \cite{s3}). La propri\'et\'e qui nous  
int\'eresse tout particuli\`erement  
est que la cat\'egorie de mod\`eles de $1$-pr\'ecat\'egories de Segal est  
une cmf \textit{interne} (voir Thm. \ref{t2}), propri\'et\'e qui n'est plus vraie  
pour la cmf des cat\'egories simpliciales. Cette propri\'et\'e nous  
permettra de d\'efinir des $Hom$ internes d\'eriv\'es entre 
$1$-cat\'egories de Segal qui seront des objets de toute importance par la suite  
(ils sont par exemple d\'ej\`a utilis\'es pour \'enoncer notre th\'eor\`eme principal). \\ 
  
\begin{df}\label{d4} 
\begin{itemize} 
\item  
Une $1$-pr\'ecat\'egorie de Segal est la donn\'ee d'un foncteur 
$$\begin{array}{cccc} 
A : & \Delta^{o} & \longrightarrow & SEns \\ 
    &   [p]      & \mapsto & A_{p} 
\end{array}$$ 
tel que $A_{0}$ soit un ensemble simplicial constant, appel\'e l'ensemble 
d'objets de $A$. 
 
\item Un morphisme entre $1$-pr\'ecatgories de Segal 
est une transformation naturelle entre foncteurs de $\Delta^{o}$ 
vers $SEns$. 
 
\item Une $1$-pr\'ecat\'egorie de Segal $A$ est une $1$-cat\'egorie de Segal si 
pour tout $[p]$ le morphisme de Segal  
$$A_{p} \longrightarrow A_{1}\times_{A_{0}}A_{1} \times \dots \times_{A_{0}} A_{1}$$ 
est une \'equivalence d'ensembles simpliciaux. 
 
\item Pour toute $1$-cat\'egorie de Segal $A$ on d\'efinit  
sa cat\'egorie homotopique $Ho(A)$ comme la cat\'egorie 
ayant $A_{0}$ comme ensemble d'objets, et $\pi_{0}(A_{1})$  
comme ensemble de morphismes. 
 
\item Un morphisme de $1$-cat\'egories de Segal $f : A \longrightarrow B$  
est une \'equivalence s'il v\'erifie les deux conditions 
suivantes. 
\begin{enumerate} 
\item Pour tout $[p] \in \Delta^{o}$,  
le morphisme $f_{p} : A_{p} \longrightarrow B_{p}$  
est une \'equivalence d'ensembles simpliciaux. 
 
\item Le morphisme induit $Ho(f) : Ho(A) \longrightarrow Ho(B)$  
est une \'equivalence de cat\'egories. 
\end{enumerate}  
\end{itemize} 
La cat\'egorie des $1$-pr\'ecat\'egories de  
Segal sera not\'ee $1-PrCat$. 
\end{df} 
 
Pour une $1$-pr\'ecat\'egorie de Segal $A$, et  
$(x,y)\in A^{2}_{0}$ une paire d'objets, on notera  
$A_{(x,y)}$ le sous-ensemble simplicial de $A_{1}$ des objets de source $x$ et but $y$  
(\cite[p. $17$]{s1}). Remarquons alors 
que l'on a  
$$A_{1}=\coprod_{(x,y)\in A_{0}^{2}}A_{(x,y)}.$$  
Il nous arrivera tr\`es souvent de d\'efinir une $1$-cat\'egorie de Segal par un  
ensemble d'objets $A_{0}$, des 
ensembles simpliciaux $A_{(x,y)}$ pour tout $(x,y)\in A_{0}^{2}$, et des  
compositions (associatives et unitaires) 
$$A_{(x,y)}\times A_{(y,z)} \longrightarrow A_{(x,z)}.$$ 
C'est \`a dire comme une cat\'egorie enrichie dans $SEns$  
(i.e. une cat\'egorie simpliciale). L'objet   
$A_{p}$ est alors d\'efini par 
$$A_{p}:=\coprod_{(x_{0}, \dots, x_{p})\in A_{0}^{p+1}} 
A_{(x_{0},x_{1})}\times A_{(x_{1},x_{2})}\times \dots 
\times A_{(x_{p-1},x_{p})}.$$ 
Remarquons que cette construction permet d'identifier la cat\'egorie des  
cat\'egories simpliciales \`a la 
sous-cat\'egorie pleine de $1-PrCat$ form\'ee des objets dont  
les morphismes de Segal sont 
des isomorphismes. Il existe donc un foncteur pleinement  
fid\`ele $S-Cat \longrightarrow 1-PrCat$,  
o\`u $S-Cat$ est la cat\'egorie des cat\'egories simpliciales.  
A l`aide de cette inclusion, on verra toute cat\'egorie simpliciale comme 
un objet de $1-PrCat$. Il est important de noter, mais nous n'utiliserons pas  
ceci par la suite, que cette inclusion induit  
une \'equivalence sur les cat\'egories homotopiques (voir par exemple  
le paragraphe \textit{Strictification} of \cite[p. $7$]{s3}). 
En particulier, ceci implique que toute $1$-cat\'egorie de Segal est  
\'equivalente \`a une cat\'egorie simpliciale. 
 
Si $M$ est une cmf simpliciale, on peut former  
la cat\'egorie simpliciale des objets cofibrants et fibrants dans $M$.  
On la consid\'erera alors comme un objet de $1-PrCat$,  
que l'on notera $E(M)$. Remarquer qu'il existe une \'equivalence  
naturelle de cat\'egories $Ho(E(M))\simeq Ho(M)$ (o\`u le membre de  
gauche est la cat\'egorie homotopique au sens de  
de la d\'efinition \ref{d4}, et celui de droite 
la cat\'egorie obtenu en localisant $M$ le long des \'equivalences).  
 
Supposons maintenant que $F : M \longrightarrow N$ soit un foncteur de Quillen  
\`a droite entre deux cmf simpliciales, pr\'eservant les 
objets cofibrants. Notons 
$G : N \longrightarrow M$ son adjoint \`a gauche, et supposons que  
$G$ soit un foncteur simplicial (i.e. de cat\'egories 
en modules sur $SEns$). Alors par adjonction on dispose  
de morphimes dans $SEns$, fonctoriels en $A$ et $B$ 
$$Hom(A,B) \longrightarrow Hom(F(A),F(B)).$$ 
Ceci permet alors de d\'efinir un foncteur simplicial de cat\'egories simpliciales 
$$E(F) : E(M) \longrightarrow E(N),$$ 
qui est consid\'er\'e comme un morphisme dans $1-PrCat$.  
 
De m\^eme, si $G : M \longrightarrow N$ est de Quillen  
\`a gauche et simplicial entre deux cmf simpliciales, et si de plus il pr\'eserve 
les objets fibrants, alors il induit un morphisme dans $1-PrCat$ 
$$E(G) : E(M) \longrightarrow E(N).$$ 
 
\begin{thm}{(\cite{s1,s5})}\label{t2} 
La cat\'egorie $1-PrCat$ est munie d'une structure  
de cmf, o\`u les cofibrations sont les monomorphismes, et les 
objets fibrants sont des $1$-cat\'egories de Segal.  
De plus, un morphisme 
entre deux $1$-cat\'egories de Segal est une \'equivalence  
dans $1-PrCat$ si et seulement si c'est une \'equivalence de $1$-cat\'egories 
de Segal (au sens de la d\'efinition \ref{d4}).  
 
Enfin, $1-PrCat$ est une cmf interne (i.e. monoidale sym\'etrique pour le produit direct, 
au sens de \cite[$\S 4.2$]{ho}).  
\end{thm} 
 
La derni\`ere assertion implique l'existence de $Hom$ internes  
$$\underline{Hom}(-,-) : 1-PrCat^{o}\times 1-PrCat \longrightarrow 1-PrCat,$$ 
satisfaisant \`a la formule d'adjonction usuelle 
$$Hom(A\times B,C)\simeq Hom(A,\underline{Hom}(B,C)),$$ 
et qui sont compatibles avec la structure de cat\'egorie de mod\`eles. En particulier, 
ils peuvent \^etre d\'eriv\'es comme il est expliqu\'e dans \cite[$\S 4.3.2$]{ho}. 
Ainsi, la cat\'egorie homotopique $Ho(1-PrCat)$ devient alors une cat\'egorie avec  
$Hom$ internes. Rappelons qu'ils sont d\'efinis 
explicitement pour $A,B \in 1-PrCat$ par la formule suivante 
$$\mathbb{R}\underline{Hom}(A,B):=\underline{Hom}(A,B') \in Ho(1-PrCat),$$ 
o\`u $B \longrightarrow B'$ est un mod\`ele fibrant,  
et qu'ils v\'erifient la formule d'adjonction 
$$\mathbb{R}\underline{Hom}(A\times B,C)\simeq \mathbb{R}\underline{Hom}(A,\mathbb{R}\underline{Hom}(B,C)).$$ 
Enfin, l'ensemble des classes d'isomorphie d'objets de  
$Ho(\mathbb{R}\underline{Hom}(A,B))$ est en bijection avec  
l'ensemble des morphismes entre $A$ et $B$ dans $Ho(1-PrCat)$ (\cite[$\S 2$]{s1}). \\ 
 
\begin{df}\label{d5} 
Soit $f : A \longrightarrow B$ un morphisme entre deux $1$-cat\'egories de Segal. 
On dit que $f$ est pleinement fid\`ele si pour tout $(x,y) \in A_{0}^{2}$ le morphisme induit 
$$f_{x,y} : A_{(x,y)} \longrightarrow B_{(f(x),f(y))}$$ 
est une \'equivalence dans $SEns$. 
 
Soit $f : A \longrightarrow B$ dans $Ho(1-PrCat)$. On dira que  
$f$ est pleinement fid\`ele si  
un repr\'esentant $f : A \longrightarrow B'$, o\`u $B \longrightarrow B'$ est un mod\`ele 
fibrant, est pleinement fid\`ele. 
\end{df} 
 
Remarquer que dans la d\'efinition pr\'ec\'edente on utilise  
le fait qu'un objet fibrant de $1-PrCat$ est automatiquement une $1$-cat\'egorie 
de Segal. \\ 
 
D'apr\`es \cite[$\S 5$]{s2} il existe pour tout $A \in Ho(1-PrCat)$ une  
famille des fl\`eches (qui est un morphisme dans 
$Ho(1-PrCat)$) 
$$Ar_{A} : A^{o}\times A \longrightarrow Top,$$ 
o\`u $Top=LSEns$ est la localisation de Dwyer-Kan de la cmf des  
ensembles simpliciaux (voir le paragraphe suivant). On peut par exemple d\'efinir 
ce morphisme en utilisant que $A$ est isomorphe dans $Ho(1-PrCat)$ 
\`a une cat\'egorie simpliciale, et en utilisant le morphisme de Yoneda enrichie pour 
les cat\'egories simpliciales. 
Par adjonction ce morphisme induit deux  
nouveaux morphismes dans $Ho(1-PrCat)$ 
$$h : A \longrightarrow \mathbb{R}\underline{Hom}(A^{o},Top)$$ 
$$k : A^{o} \longrightarrow \mathbb{R}\underline{Hom}(A,Top).$$ 
 
Le lemme de Yoneda pour les $1$-cat\'egories de Segal est le  
th\'eor\`eme fondamental suivant. 
 
\begin{thm}{(\cite[Thm. $2$]{s3})}\label{t3} 
Pour tout $A \in Ho(1-PrCat)$ les morphismes $h$ et $k$ sont pleinement fid\`eles.  
\end{thm} 
 
Remarquons que la notion de $1$-pr\'ecat\'egorie de  
Segal est une g\'en\'eralisation de celle de pr\'e-$\Delta^{o}$-espace  
(qui correspond au cas 
o\`u $A_{0}=*$). Plus pr\'ecis\'ement, il existe un foncteur \'evident 
$$B : Pr\Delta^{o}-SEns \longrightarrow */1-PrCat,$$ 
qui identifie $Pr\Delta^{o}-SEns$ avec la sous-cat\'egorie pleine  
form\'ee des objets $A$ tels que $A_{0}=*$. Il est facile de voir  
qu'il poss\`ede un adjoint \`a droite 
$$\Omega_{*} : */1-PrCat \longrightarrow Pr\Delta^{o}-SEns.$$ 
Pour un objet point\'e $(A,a) \in */1-PrCat$ il est  
d\'efini explicitement de la fa\c{c}on suivante 
$$\begin{array}{cccc} 
\Omega_{a}(A) : & \Delta^{o} & \longrightarrow & SEns \\ 
& [m] & \mapsto & A_{\underbrace{(a,a,\dots,a)}_{m+1 \; fois}}. 
\end{array}$$ 
D'apr\`es \cite[$\S 3.11$]{s5} il est \'evident que le  
foncteur $B$ est de Quillen \`a gauche, et donc que $\Omega_{*}$ est de Quillen \`a droite. On en d\'eduit 
donc un foncteur d\'eriv\'e 
$$\mathbb{R}\Omega_{*} : Ho(*/1-PrCat) \longrightarrow Ho(Pr\Delta^{o}-SEns).$$ 
En utilisant le fait qu'un objet fibrant de $1-PrCat$ est  
automatiquement une $1$-cat\'egorie de Segal, on remarque que  
$\mathbb{R}\Omega_{*}$ prend ses valeurs dans la sous-cat\'egorie pleine des $\Delta^{o}$-espaces. 
 
\begin{df}\label{d6} 
Soit $A \in Ho(1-PrCat)$, et $f$ un objet de $\mathbb{R}\underline{Hom}(A,Top)$.  
On d\'efinit le $\Delta^{o}$-espace des  
endomorphismes de $f$ par 
$$\mathbb{R}\underline{End}(f):=\mathbb{R}\Omega_{f}(\mathbb{R}\underline{Hom}(A,Top)) \in Ho(Pr\Delta^{o}-SEns).$$ 
\end{df} 
 
\end{subsection} 
 
\begin{subsection}{Localisation de Dwyer-Kan} 
 
Soit $C$ une cat\'egorie, et $S$ un ensemble de morphismes.  
D'apr\`es \cite{dk} il existe une cat\'egorie simpliciale $L(C,S)$, qui  
d\'epend fonctoriellement du couple $(C,S)$, et telle que  
la cat\'egorie homotopique de $L(C,S)$ soit \'equivalente \`a la cat\'egorie 
localis\'ee $S^{-1}C$. Pour tout foncteur  
$F : (C,S) \longrightarrow (D,T)$ (i.e. un foncteur $F : C \longrightarrow D$ 
qui envoit $S$ dans $T$) nous noterons 
$LF : L(C,S) \longrightarrow L(D,T)$ le morphisme  
induit dans $1-PrCat$. D'un point de vue ensembliste remarquons que si  
$\mathbb{U}$ est un univers et $C$ une 
cat\'egorie avec $C\in \mathbb{U}$,  
alors pour tout ensemble $S$ de morphismes de $C$,  
$L(C,S)$ appartient encore \`a $\mathbb{U}$. 
 
Lorsque $C$ est une cmf et que $S$ est l'ensemble des  
\'equivalences nous \'ecrirons simplement $LC$ pour $L(C,S)$. De m\^eme, pour toute 
sous-cat\'egorie pleine $C$ d'une cmf $M$, nous \'ecrirons  
$LC$ pour $L(C,Equiv.\cap C)$, o\`u $Equiv.$ est l'ensemble des \'equivalences de 
$M$. 
 
Les trois propri\'et\'es fondamentales que nous utiliserons implicitement sont les suivantes. 
 
\begin{itemize} 
\item 
Si $F : (C,S) \longrightarrow (D,T)$ et  
$G : (D,T) \longrightarrow (C,S)$ sont deux foncteurs tels qu'il  
existe des diagrammes de transformations naturelles  
$$\xymatrix{Id \ar[r] & A & \ar[l] FG} \qquad \xymatrix{GF & \ar[l] B \ar[r] & Id}$$ 
qui sont objet par objet des morphismes de $T$ et  
de $S$, alors $LF$ et $LG$ sont inverses l'un de l'autre dans $Ho(1-PrCat)$ (voir  
\cite[$\S 8.1$]{s1}). En particulier, si $M^{c}$  
est la cat\'egorie des objets cofibrants dans une cmf $M$, alors le morphisme naturel 
$$LM^{c} \longrightarrow  LM,$$ 
est un isomorphisme dans $Ho(1-PrCat)$. En effet, le  
foncteur de remplacement cofibrant $Q : M \longrightarrow M^{c}$ 
induit un morphisme inverse dans $Ho(1-PrCat)$,  
$LQ : LM \longrightarrow LM^{c}$ (\cite[$\S 1.1$]{ho}).  
De m\^eme, si $M^{f}$ est la cat\'egorie des objets  
fibrants dans $M$, alors le morphisme naturel  
$$LM^{f} \longrightarrow LM$$ 
est un isomorphisme dans $Ho(1-PrCat)$. 
Ou encore, si $M^{cf}$ est la sous-cat\'egorie pleine  
des objets fibrants et cofibrants de $M$, le morphisme naturel 
$$LM^{cf} \longrightarrow LM,$$ 
est un isomorphisme dans $Ho(1-PrCat)$. 
 
\item 
Soit $M$ une cmf simpliciale, et $E(M)$ la cat\'egorie  
simpliciale des objets fibrants et cofibrants de $M$. Alors il existe un 
isomorphisme naturel dans $Ho(1-PrCat)$ 
$$E(M) \simeq LM.$$ 
  
Soit $G : M \longrightarrow N$ un foncteur de Quillen  
\`a gauche et simplicial entre deux cmf simpliciales, pr\'eservant les 
objets fibrants et les \'equivalences. Alors il existe  
un isomorphisme dans $Ho(*/1-PrCat)$ 
$$(\mathbb{R}\underline{Hom}(LM,LN),LG)\simeq (\mathbb{R}\underline{Hom}(E(M),E(N)),E(G)).$$ 
 
De m\^eme, si $F : M \longrightarrow N$ est de Quillen  
\`a droite avec un adjoint \`a gauche simplicial, et s'il pr\'eserve les 
objets cofibrants et les \'equivalences, alors il existe  
un isomorphisme dans $Ho(*/1-PrCat)$ 
$$(\mathbb{R}\underline{Hom}(LM,LN),LF)\simeq (\mathbb{R}\underline{Hom}(E(M),E(N)),E(F)).$$ 
 
\item Si $M_{0}$ est une sous-cat\'egorie pleine d'une  
cmf $M$, qui est stable par \'equivalence (i.e. tout objet de $M$ isomorphe dans $Ho(M)$ 
\`a un objet de $M_{0}$ est dans $M_{0}$), alors le  
morphisme $LM_{0} \longrightarrow LN$ est pleinement fid\`ele 
(au sens de la d\'efinition \ref{d5}). 
 
\end{itemize} 
 
Remarquons que la seconde propri\'et\'e permet d'associer  
\`a tout foncteur de Quillen \`a droite entre deux cmf $F : M \longrightarrow N$,  
un morphisme $LF : LM^{f} \longrightarrow LN^{f}$, ou encore un objet  
$$(\mathbb{R}\underline{Hom}(LM^{f},LN^{f}),LF) \in Ho(*/1-PrCat).$$ 
En effet, d'apr\`es \cite[$\S 1.1.12$]{ho} le foncteur  
$F : M^{f} \longrightarrow N^{f}$ pr\'eserve les \'equivalences.  
Mais comme $LM^{f}\simeq LM$ et $LN^{f}\simeq LN$ ceci  
d\'etermine un objet bien d\'efini \`a isomorphisme unique pr\`es 
$$(\mathbb{R}\underline{Hom}(LM,LN),LF) \in Ho(*/1-PrCat).$$ 
De m\^eme, si $G : M \longrightarrow N$ est de Quillen \`a gauche on en d\'eduit un objet  
$$(\mathbb{R}\underline{Hom}(LM,LN),LG) \in Ho(*/1-PrCat).$$ 
Il est important de remarquer qu'un foncteur de Quillen  
\`a droite (resp. \`a gauche) $F : M \longrightarrow N$ qui est une  
\'equivalence de Quillen induit un isomorphisme dans $Ho(1-PrCat)$,  
$LF : LM \longrightarrow LN$. De plus, si $G$ est l'adjoint 
\`a gauche (resp. \`a droite) de $F$, alors $LG$ est l'inverse de $LF$. \\ 
 
Terminons par l'\'etude d'un cas particulier dans lequel  
on peut calculer le $\Delta^{o}$-espace des endomorphismes d'un  
morphisme $A \longrightarrow Top$ \`a l'aide du lemme de 
Yoneda. Cet exemple nous sera tr\`es utile lors de la preuve du th\'eor\`eme principal.  
 
Supposons que $F : M \longrightarrow SEns$ soit un foncteur  
de Quillen \`a droite, o\`u $M$ est une cmf simpliciale, et que 
son adjoint \`a gauche $G : SEns \longrightarrow M$ pr\'eserve  
la structure simpliciale. On a vu que l'on pouvait construire un morphisme 
dans $1-PrCat$ 
$$E(F) : E(M) \longrightarrow E(SEns).$$ 
Supposons que $G$ pr\'eserve les objets fibrants, et notons  
$X=G(*) \in E(M)$. Par adjonction, on a $F(Y)=\underline{Hom}(X,Y)$, o\`u  
$\underline{Hom}$ est le $Hom$ simplicial de $M$. A l'aide de la description 
de $Ar_{E(M)}$ (\cite[$\S 5.6$]{s2}) cela implique facilement que  
l'objet $X$ est envoy\'e par $k$ sur l'image de $E(F)$ dans  
$\mathbb{R}\underline{Hom}(A,E(SEns))$. 
On obtient donc un morphisme dans $Ho(*/1-PrCat)$ 
$$k : (E(M)^{o},X) \longrightarrow (\mathbb{R}\underline{Hom}(E(M),E(SEns)),E(F)).$$ 
D'apr\`es le lemme de Yoneda \ref{t3} et la d\'efinition  
de pleinement fid\`ele, on en d\'eduit que le morphisme induit dans $Ho(Pr\Delta^{o}-SEns)$ 
$$\mathbb{R}\Omega_{X}(E(M))^{o}\longrightarrow \mathbb{R}\Omega_{E(F)}(\mathbb{R}\underline{Hom}(E(M),E(SEns)))$$ 
est un isomorphisme. Enfin, il est clair que  
$\mathbb{R}\Omega_{X}(E(M))^{o}$ est isomorphe \`a  
$\underline{Hom}(X,X)^{o}$, le mono\"{\i}de simplicial oppos\'e  
des endomorphismes de $X$ dans $M$. En appliquant  
les propri\'et\'es de la localisation de Dwyer-Kan rappel\'ees ci-dessus on obtient le  
corollaire suivant. 
 
\begin{cor}\label{c1} 
Soit $F : M \longrightarrow SEns$ un foncteur de Quillen  
\`a droite pr\'eservant les \'equivalences.  
On suppose que son adjoint \`a gauche $G$ pr\'eserve la  
structure simpliciale et les objets fibrants. Soit $X=G(*) \in M$. Alors, 
il existe un isomorphisme naturel dans $Ho(Pr\Delta^{o}-SEns)$ 
$$\underline{Hom}(X,X)^{o} \simeq \mathbb{R}\underline{End}(LF).$$ 
\end{cor} 
 
\end{subsection} 
 
\end{section} 
 
\begin{section}{Classifications des champs localement constants} 
 
Ce chapitre est essentiellement compos\'e de trois parties.  
Dans un premier temps nous d\'efinirons la $1$-cat\'egorie de Segal $Loc(X)$, des 
champs localement constants sur un $CW$ complexe $X$, ainsi  
que son foncteur fibre en $x \in X$, $\omega_{x}$. Nous introduirons  
ensuite la $1$-cat\'egorie de Segal $Fib(S)$, des fibrations sur un ensemble  
simplicial $S$. Par des techniques d'alg\`ebre homotopique nous montrerons 
que $Loc(X)$ est naturellement \'equivalente \`a $Fib(S(X))$,  
o\`u $S(X)$ est l'ensemble simplicial singulier de $X$ (voir  
le th\'eor\`eme \ref{p3}).  
De plus, cette \'equivalence est compatible aux foncteurs fibres.  
C'est ce que nous r\'esumons en affirmant que \textit{les champs localement  
constants sur $X$ sont classifi\'es par les fibrations sur  
$S(X)$}. Ceci suppose connue l'existence 
d'une structure de cat\'egorie de mod\`eles ferm\'ee sur les  
pr\'efaisceaux simpliciaux sur $X$ pour laquelle les \'equivalences sont 
les \'equivalences fibre \`a fibre. Nous avons choisi de travailler  
avec la structure projective de \cite[$\S 5$]{s1} et \cite{bl} car ceci permet 
d'utiliser qu'une certaine adjonction est de Quillen (voir la  
proposition \ref{p2}), ce qui n'est plus vrai pour la structure d\'efinie  
par J. F. Jardine dans \cite{j}. Enfin, nous d\'efinirrons pour tout groupe  
simplicial $G$ la $1$-cat\'egorie de Segal de ses repr\'esentations, $G-Top$.  
Le second th\'eor\`eme de classification affirme alors que  
la $1$-cat\'egorie de Segal $G-Top$ est naturellement \'equivalente \`a 
$Fib(d(G))$, o\`u $d(G)$ est l'ensemble simplicial classifiant  
du groupe simplicial $G$ (voir le th\'eor\`eme \ref{p5}).  
 
\begin{subsection}{Champs localement constants} 
 
Fixons $\mathbb{U}$ un univers avec $\Delta \in \mathbb{U}$,  
et $\mathbb{V}$ un univers avec  
$\mathbb{U} \in \mathbb{V}$. Nous noterons $Es^{cw}_{\mathbb{U}}$  
la sous-cat\'egorie pleine de 
$Es_{\mathbb{U}}$ form\'ee des $CW$ complexes. Remarquons que  
les objets de $Es_{\mathbb{U}}^{cw}$ sont tous localement contractiles et 
h\'er\'editairement paracompacts (voir les th\'eor\`emes  
$1.3.2$, $1.3.5$ et Ex. $1$ du $\S 1.3$ de \cite{fp}).  
 
Pour $X \in Es_{\mathbb{U}}$, on note $SPr(X)$ la cat\'egorie des pr\'efaisceaux sur 
$X$ \`a valeurs dans $SEns_{\mathbb{U}}$.  
On munit $SPr(X)$ de la structure de cmf de \cite[$\S 5$]{s1} 
(voir aussi \cite{bl}).  
Nous mettons en garde le lecteur qu'il ne s'agit pas 
de la structure d\'efinie dans \cite{j}, mais de la structure de type $HBKQ$ de \cite[$\S 5$]{s1},  
o\`u les cofibrations sont engendr\'ees par \textit{les  
additions libres de cellules}. Cependant les \'equivalences dans $SPr(X)$ sont aussi les  
\'equivalences utilis\'ees dans \cite{j}. Dans notre  
situtation ce sont donc exactement les morphismes $f : F \longrightarrow G$ tels que  
pour tout $x \in X$, le morphisme induit sur les  
fibres $f_{x} : F_{x} \longrightarrow G_{x}$ soit une \'equivalence.  
 
Rappelons aussi que la structure de cmf sur $SPr(X)$ est obtenue 
par localisation de Bousfield \`a gauche de la structure de cmf pour  
la topologie grossi\`ere sur la cat\'egorie $Ouv(X)$ 
des ouverts de $X$ en \textit{inversant les hyper-recouvrements}.  
Une d\'emonstration de ce dernier r\'esultat pourra se 
trouver dans le travail \cite{du}. Ceci entraine en  particulier une jolie charact\'erisation 
des objets fibrants. Afin de l'\'enoncer rappelons la d\'efinition d'hyper-recouvrement 
de \cite[$8.4$]{am}. 
 
\begin{df}\label{d6'} 
Un ($\mathbb{U}$-)hyper-recouvrement $U_{*}$ d'un ouvert $U \subset X$ est  
un faisceau simplicial $U_{*}$ sur l'espace topologique $U$ satisfaisant aux 
propri\'et\'es suivantes. 
\begin{enumerate} 
\item Pour tout entier $n\geq 0$, le faisceau $U_{n}$ est isomorphe \`a une 
r\'eunion disjointe d'une famille $\mathbb{U}$-petite de faisceaux  
repr\'esentables par des ouverts de $U$. On \'ecrira symboliquement 
$$U_{n}\simeq \coprod_{i \in I_{n}} U_{n}^{(i)},$$ 
avec $I_{n} \in \mathbb{U}$ et $U_{n}^{(i)}$ des ouverts de $U$. 
 
\item Le morphisme $U_{0} \longrightarrow *$ est un 
\'epimorphisme de faisceaux sur $U$. En d'autres termes 
la famille $\{U_{0}^{(i)}\}_{i \in I_{0}}$ est un 
recouvrement ouvert de $U$. 
 
\item Pour tout entier $n\geq 0$, le morphisme 
$$U_{n+1} \longrightarrow (Cosq_{n}U_{*})_{n+1}$$ 
est un \'epimorphisme de faisceaux sur $U$. 
 
\end{enumerate} 
\end{df} 
 
Dans la d\'efinition pr\'ec\'edente, $(Cosq_{n}U_{*})_{n+1}$ d\'esigne  
le $(n+1)$-\'etage du $n$-\`eme co-squelette 
de l'objet simplicial $U_{*}$, qui est d\'efini dans \cite[$\S 8$]{am}, 
et que l'on peut aussi d\'efinir comme suit. Comme la cat\'egorie 
des faisceaux sur $U$ (\`a valeurs dans $\mathbb{U}$) est compl\'ete et co-compl\`ete, on peut munir 
la cat\'egorie des faisceaux simpliciaux sur $U$ d'une structure de cat\'egorie 
enrichie sur $SEns_{\mathbb{U}}$. De plus, cette structure enrichie poss\`ede des produits externes 
et des objets exponentiels par des objets de $SEns_{\mathbb{U}}$ (voir par exemple \cite[$\S II.2$]{gj}). 
En particulier, pour un faisceau simplicial $U_{*}$, le faisceau simplicial $U_{*}^{\partial \Delta^{n+1}}$ 
existe. On a alors 
$$(Cosq_{n}U_{*})_{n+1}=(U_{*}^{\partial \Delta^{n+1}})_{0}.$$ 
 
Dans le lemme suivant, on utilisera la notation $F(U_{n})$ pour $F$ un pr\'efaisceau simplicial sur $X$ et 
$U_{*}$ un hyper-recouvrement d'un ouvert $U$. Avec les notations de la d\'efinition \ref{d6'}, on 
d\'esigne de cette fa\c{c}on 
$$F(U_{n}):=\prod_{i \in I_{n}}F(U_{n}^{(i)}).$$ 
 
\begin{lem}\label{l4'}{(Voir \cite{du})} 
Un objet $F \in SPr(X)$ est fibrant si et seulement si les deux assertions suivantes sont satisfaites. 
 
\begin{enumerate} 
\item Pour tout ouvert $U \subset X$, $F(U)$ est fibrant dans $SEns_{\mathbb{U}}$.  
 
\item Pour tout ouvert $U \subset X$, et tout hyper-recouvrement ouvert $U_{*}$ de $U$, le morphisme naturel 
$$F(U) \longrightarrow Holim_{[m] \in \Delta}F(U_{m})$$ 
est une \'equivalence dans $SEns_{\mathbb{U}}$. 
 
\end{enumerate} 
\end{lem} 
 
\textit{Remarque:} Un pr\'efaisceau simplicial qui v\'erifie la  
seconde condition est g\'en\'eralement appel\'e  
un \textit{champ} (c'est le point de vue utilis\'e dans \cite{s1}).  
De plus, bien que la th\'eorie ne soit pas encore  
d\'evelopp\'ee, la th\'eorie homotopique des pr\'efaisceaux simpliciaux 
est sens\'ee \^etre \'equivalente \`a celle des pr\'efaisceaux  
en $\infty$-groupo\"{\i}des. Ainsi, comme la $1$-cat\'egorie de Segal $LSPr(X)$ est  
\'equivalente \`a la cat\'egorie simpliciale des objets  
fibrants et cofibrants dans $SPr(X)$, on peut l\'egitimement penser \`a  
$LSPr(X)$ comme \`a un mod\`ele pour la $\infty$-cat\'egorie  
des $\infty$-champs en groupo\"{\i}des sur $X$. Nous appellerons alors 
\textit{champ} tout objet de $LSPr(X)$.  
 
\begin{df}\label{d6''} 
La $1$-cat\'egorie de Segal des champs sur $X$ est $LSPr(X)$. 
\end{df} 
 
Rappelons aussi que $SPr(X)$ est aussi une cmf simpliciale,  
le produit externe de $F \in SPr(X)$ par $A \in SEns_{\mathbb{U}}$ \'etant d\'efini 
par la formule  
$$\begin{array}{cccc} 
A\times F : & Ouv(X)^{o} & \longrightarrow & SEns_{\mathbb{U}} \\ 
& U & \mapsto & A\times F(U). 
\end{array}$$ 
 
On dispose alors de $Hom$ \`a valeurs dans $SEns_{\mathbb{U}}$, qui seront 
not\'es comme d'habitude $\underline{Hom}$. Par d\'efinition, on a  
$\underline{Hom}(F,G)_{m}:=Hom(\Delta^{m}\times F,G)$. Enfin, il existe une 
op\'eration d'exponentiation, qui \`a $F \in SPr(X)$ et  
$A \in SEns_{\mathbb{U}}$ associe 
le pr\'efaisceau simplicial $F^{A}$ d\'efini par  
$$\begin{array}{cccc} 
F^{A} : & Ouv(X)^{o} & \longrightarrow & SEns_{\mathbb{U}} \\ 
& U & \mapsto & F(U)^{A}. 
\end{array}$$ 
Ces constructions satisfont aux relations  
d'adjonctions usuelles, et de plus aux axiomes de \cite[$\S 4.2.18$]{ho} d\'efinissant  
la structure de cmf simpliciale. \\ 
 
Pour tout morphisme dans $Es_{\mathbb{U}}^{cw}$, $f : Y \longrightarrow X$,  
on dispose d'un foncteur d'images directes 
$$f_{*} : SPr(Y) \longrightarrow SPr(X)$$ 
d\'efini par la formule $f_{*}(F)(U)=F(f^{-1}(U))$.  
Ce foncteur poss\`ede un adjoint \`a gauche 
$$f^{*} : SPr(X) \longrightarrow SPr(Y).$$ 
Comme $f^{*}$ pr\'eserve les \'equivalences et les cofibrations,  
$f^{*}$ est de Quillen \`a gauche, et donc $f_{*}$ est de Quillen \`a droite.  
 
Notons que $X \mapsto SPr(X)$, $f \mapsto f^{*}$  
d\'efinit une cat\'egorie cofibr\'ee sur $Es_{\mathbb{U}}$. En appliquant le proc\'ed\'e standard 
(voir par exemple \cite[Thm. $3.4$]{ma})
de strictification on supposera que ceci d\'efinit un vrai foncteur 
$$\begin{array}{ccc} 
Es_{\mathbb{U}} & \longrightarrow & Cat_{\mathbb{V}} \\ 
X & \mapsto & SPr(X) \\ 
(f : X \rightarrow Y) & \mapsto & (f^{*} : SPr(Y) \rightarrow SPr(X)) 
\end{array}$$ 
On se permettra donc de supposer par la suite que  
pour $f  : Y \longrightarrow X$ et $g : Z \longrightarrow Y$ deux morphismes, on a 
une \'egalit\'e $g^{*}f^{*}=(fg)^{*}$. \\ 
 
\textit{Remarque:} Le choix des univers $\mathbb{U}$ et  
$\mathbb{V}$ est tel qu'apr\`es strictification les  
cat\'egories $SPr(X)$ appartiennent toujours 
\`a $\mathbb{V}$. De plus, comme la notion de structure  
de cmf est invariante par \'equivalence de cat\'egories, les cat\'egories 
$SPr(X)$ sont munies naturellement de structures de cmf.  
De m\^eme, les foncteurs $f^{*}$ restent de Quillen \`a gauche. 
 
\begin{df}\label{d7} 
Soit $X \in Es_{\mathbb{U}}$. 
\begin{itemize} 
 
\item Un objet $F \in SPr(X)$ est $h$-constant,  
s'il est isomorphe dans $Ho(SPr(X))$ \`a un pr\'efaisceau constant.  
 
\item Un objet $F \in SPr(X)$ est localement  
$h$-constant s'il existe un recouvrement ouvert $\{U_{i}\}_{i \in I}$ de $X$, tel que 
pour tout $i \in I$ la restriction de $F$ \`a  
$U_{i}$, $F_{|U_{i}} \in SPr(U_{i})$ soit $h$-constant.  
 
\item La sous-cat\'egorie pleine de $SPr(X)$  
form\'ee des objets localement $h$-constants sera not\'ee  
$PrLoc(X)$.  
 
\item On d\'efinit 
$$Loc(X):=LPrLoc(X).$$  
Elle sera consid\'er\'ee comme un objet de  
$1-PrCat_{\mathbb{V}}$. 
\end{itemize} 
La $1$-cat\'egorie de Segal $Loc(X)$ est  
appel\'ee la $1$-cat\'egorie de Segal des champs localement constants sur $X$. 
\end{df} 
 
\textit{Remarque:} Comme la sous-cat\'egorie  
$PrLoc(X)$ de $SPr(X)$ est stable par \'equivalences, le foncteur naturel 
$Loc(X) \longrightarrow LSPr(X)$ est pleinement  
fid\`ele. Nous utiliserons en particulier que le foncteur induit $Ho(Loc(X)) \longrightarrow 
Ho(SPr(X))$ est pleinement fid\`ele, et que son  
image essentielle est form\'ee des objets localement $h$-constants. \\ 
 
Rappelons que pour $X \in Es_{\mathbb{U}}$,  
$F \in SPr(X)$ et $x \in X$, on d\'efinit comme  
il en est l'usage la fibre en $x$ par 
$$F_{x}:=Colim_{x \in U}F(U)$$ 
o\`u la colimite est prise sur les voisinages de $x$ dans $X$.  
 
On dispose alors du foncteur fibre en $x$ 
$$\omega_{x} : SPr(X) \longrightarrow SEns_{\mathbb{U}}.$$ 
Comme c'est un foncteur qui pr\'eserve les  
\'equivalences il induit un morphisme dans $1-PrCat_{\mathbb{V}}$ 
$$L\omega_{x} : LPrLoc(X)=Loc(X) \longrightarrow LSEns_{\mathbb{U}}=Top.$$ 
On obtient de cette fa\c{c}on un objet  
$$(\mathbb{R}\underline{Hom}(Loc(X),Top),L\omega_{x}) \in Ho((1-PrCat_{\mathbb{V}})_{*}).$$ 
 
Soit $f : (Y,y) \longrightarrow (X,x)$  
un morphisme dans $*/Es_{\mathbb{U}}$. Le foncteur d'image inverse  
$f^{*}$ pr\'eserve clairement les  
\'equivalences ainsi que les objets localement $h$-constants.  
Il induit donc un morphisme sur les localis\'ees de Dwyer-Kan 
$$Lf^{*} : LPrLoc(X)=Loc(X) \longrightarrow LPrLoc(Y)=Loc(Y).$$ 
Ce morphisme \'etant compatible avec les foncteurs  
fibres en $x$ et $y$ il induit un morphisme dans $Ho(*/1-PrCat_{\mathbb{V}})$ 
$$(\mathbb{R}\underline{Hom}(Loc(Y),Top),L\omega_{y}) \longrightarrow  
(\mathbb{R}\underline{Hom}(Loc(X),Top),L\omega_{x}).$$ 
On d\'efinit ainsi un foncteur 
$$\begin{array}{cccc}  
(Loc(-),L\omega) : & */Es^{cw}_{\mathbb{U}} & \longrightarrow & Ho(*/1-PrCat_{\mathbb{V}}) \\ 
& (X,x) & \mapsto & (\mathbb{R}\underline{Hom}(Loc(X),Top),L\omega_{x}). 
\end{array}$$ 
Il est important de remarquer que nous n'avons d\'efini  
le foncteur ci-dessus que sur la sous-cat\'egorie de $Es$ form\'ee des $CW$ complexes.  
 
\begin{lem}\label{l3'} 
Soit $F$ un objet dans $PrLoc(X)$, et  
$f : Y \longrightarrow X$ un morphisme dans $Es_{\mathbb{U}}$.  
Alors pour tout $y \in Y$ le morphisme  
induit sur les fibres 
$$f^{*}(F)_{y} \longrightarrow F_{f(y)}$$ 
est un isomorphisme. 
\end{lem} 
 
\textit{Preuve:} C'est imm\'ediat par d\'efinition. \hfill $\Box$ \\ 
 
\begin{cor}\label{c3'} 
Soit $f : Y \longrightarrow X$ un morphisme  
surjectif dans $Es_{\mathbb{U}}$, et $u : F \longrightarrow G$ un morphisme dans $Ho(SPr(X))$.  
Alors $u$ est un isomorphisme si et seulement si le morphisme induit 
$$f^{*}(F) \longrightarrow f^{*}(G)$$ 
est un isomorphisme dans $Ho(SPr(Y))$. 
\end{cor} 
 
\textit{Preuve:} En effet, comme $f$ est surjective  
le lemme \ref{l3'} entraine que pour tout $x \in X$ le morphisme 
$u_{x} : F_{x} \longrightarrow G_{x}$ est une \'equivalence. \hfill $\Box$ \\ 
 
\end{subsection} 
 
\begin{subsection}{Champs localement constants et fibrations} 
 
Dans ce paragraphe nous allons d\'emontrer un analogue  pour les 
champs localement constants de la correspondance entre faisceaux localement  
constants et rev\^etements. Pour cela nous construirons pour tout $CW$ complexe $X$,  
une adjonction de Quillen 
($Se$ pour \textit{section}, et $R$ pour \textit{r\'ealisation}) 
$$Se : SEns/S(X) \longrightarrow SPr(X) \qquad R : SPr(X) \longrightarrow SEns/S(X),$$ 
et nous montrerons que les morphismes induits sur les localis\'ees de Dwyer-Kan 
$$LSe : LSEns/S(X) \longrightarrow LSPr(X) \qquad LR : LSPr(X) \longrightarrow LSEns/S(X),$$ 
induisent des isomorphismes inverses l'un de l'autre dans $Ho(1-PrCat)$  
entre $Loc(X)$ et $LSEns/S(X)$ (voir le th\'eor\`eme \ref{p3}).  
Remarquons que $LSEns/S(X)$ est naturellement isomorphe dans  
$Ho(1-PrCat)$ \`a la 
cat\'egorie simpliciale des fibrations sur $S(X)$. Le r\'esultat que nous  
d\'emontrerons  peut donc se r\'esumer en affirmant  
que \textit{les champs localement constants sur $X$ sont  
classifi\'es par les fibrations sur $S(X)$}. \\ 
 
Soit $S\in SEns_{\mathbb{U}}$ un ensemble simplicial et $SEns_{\mathbb{U}}/S$ la cmf des ensembles simpliciaux sur $S$.  
 
\begin{df} 
Pour tout ensemble simplicial $S$, nous d\'efinissons la $1$-cat\'egorie de Segal des fibrations sur $S$ par 
$$Fib(S):=L(SEns_{\mathbb{U}}/S) \in 1-PrCat_{\mathbb{V}}.$$ 
Nous noterons aussi  
$$Top:=Fib(*)\in 1-PrCat_{\mathbb{V}}.$$ 
\end{df} 
 
Remarquer que $Fib(S)$ est une $1$-cat\'egorie de Segal \'equivalente  
\`a la cat\'egorie simpliciale des fibrations sur $S$. Ceci justifie 
le choix de la terminologie. De m\^eme $Top$ est \'equivalente \`a  
la cat\'egorie simpliciale des ensemble simpliciaux fibrants. Elle 
est donc \'equivalente \`a la localis\'ee de Dwyer-Kan des espaces  
topologiques, d'o\`u le choix de la notation $Top$. \\ 
 
Soit $f : S' \longrightarrow S$ un morphisme dans $SEns_{\mathbb{U}}$.  
On dispose alors de deux foncteurs adjoints ($f^{*}$ \'etant adjoint \`a droite
de $f_{!}$)
$$\begin{array}{cccc} 
f^{*} : & SEns_{\mathbb{U}}/S & \longrightarrow & SEns_{\mathbb{U}}/S' \\
& (Y\rightarrow S) & \mapsto & (Y\times_{S}S'\rightarrow S'), 
\end{array}$$ 
$$\begin{array}{cccc} 
f_{!} : & SEns_{\mathbb{U}}/S' & \longrightarrow & SEns_{\mathbb{U}}/S \\   
& (Z\rightarrow S') & \mapsto & (Z\rightarrow S'\rightarrow S). 
\end{array}$$
La correspondance $S \mapsto SEns_{\mathbb{U}}/S$, $f \mapsto f^{*}$  
d\'efinit une cat\'egorie cofibr\'ee sur $SEns_{\mathbb{U}}$  
\`a valeurs dans $Cat_{\mathbb{V}}$.  
Pour \^etre tout \`a fait rigoureux nous supposerons par  
la suite que nous l'avons strictifi\'ee  
et donc que pour deux morphismes 
$f : S' \longrightarrow S$ et $g : S'' \longrightarrow S'$  
on a $g^{*}f^{*}=(gf)^{*}$. \\ 
 
\textit{Remarque:} Comme dans le paragraphe pr\'ec\'edent,  
les cat\'egories $SEns_{\mathbb{U}}/S(X)$ restent des cmf de  
$\mathbb{V}$, et les foncteurs 
$f^{*}$ restent de Quillen \`a gauche. \\ 
 
Le foncteur $f^{*} : SEns_{\mathbb{U}}/S \longrightarrow SEns_{\mathbb{U}}/S'$  
poss\`ede aussi un adjoint \`a droite,  
$$f_{*} : SEns_{\mathbb{U}}/S' \longrightarrow SEns_{\mathbb{U}}/S.$$ 
En effet, $f^{*}$ commute avec les $\mathbb{U}$-colimites, et  
$SEns_{\mathbb{U}}/S'$ poss\`ede un ensemble d'objets g\'en\'erateurs 
appartenant \`a $\mathbb{U}$ (\`a savoir les morphismes $\Delta^{m} \rightarrow S'$).  
 
Notons enfin que les foncteurs $f^{*}$ et $f_{!}$ pr\'eservent les structures simpliciales,  
et donc les isomorphismes d'adjonctions induisent des isomorphismes naturels 
$$\underline{Hom}_{S}(A,f_{*}B)\simeq \underline{Hom}_{S'}(f^{*}A,B) \qquad  
\underline{Hom}_{S}(f_{!}A,B)\simeq \underline{Hom}_{S'}(A,f^{*}B),$$ 
o\`u $\underline{Hom}_{S}$ d\'esigne le $Hom$ simplicial de $SEns_{\mathbb{U}}/S$. \\  
 
Le foncteur $f_{!}$ est clairement de Quillen \`a gauche,  
et donc $f^{*}$ est de Quillen \`a droite. Ainsi, on en  
d\'eduit un morphisme dans $Ho(1-PrCat_{\mathbb{V}})$ 
$$Lf^{*} : L(SEns_{\mathbb{U}}/S)=Fib(S) \longrightarrow (LSEns_{\mathbb{U}}/S')=Fib(S').$$ 
Par exemple, lorsque $S'=\{s\}$ est un point de $S$  
on obtient un foncteur fibre 
$$L\omega_{s} : Fib(S) \longrightarrow Fib(*)=Top.$$ 
Ainsi, si $(S,s) \in */SEns_{\mathbb{U}}$ est un ensemble simplicial point\'e on en d\'eduit un objet 
$$(\mathbb{R}\underline{Hom}(Fib(S),Top),L\omega_{s}) \in Ho(*/1-PrCat_{\mathbb{V}})$$ 
qui d\'epend fonctoriellement du couple $(S,s)$. Ceci d\'etermine donc un foncteur 
$$\begin{array}{cccc}  
(Fib(-),L\omega) : & */SEns_{\mathbb{U}} & \longrightarrow & Ho(*/1-PrCat_{\mathbb{V}}) \\ 
& (S,s) & \mapsto & (\mathbb{R}\underline{Hom}(Fib(S),Top),L\omega_{s}) 
\end{array}$$ 
Remarquons enfin, que si $f$ est une fibration,  
alors $f^{*}$ pr\'eserve les \'equivalences  
(car la cmf $SEns_{\mathbb{U}}$ est propre) 
ainsi que les cofibrations et donc est de Quillen \`a gauche. 
Ceci implique par adjonction que $f_{*}$  
est alors de Quillen \`a droite (ceci n'est plus vrai pour $f$ quelconque).  
 
\begin{lem}\label{l1} 
Soit $f : (S,s) \longrightarrow (S',s')$ une \'equivalence dans $*/SEns_{\mathbb{U}}$. Alors le 
morphisme induit  
$$(\mathbb{R}\underline{Hom}(Fib(S),Top),L\omega_{s}) \longrightarrow  
(\mathbb{R}\underline{Hom}(Fib(S'),Top),L\omega_{s'})$$ 
est un isomorphisme dans $Ho(*/1-PrCat_{\mathbb{V}})$ 
\end{lem} 
 
\textit{Preuve:} Il suffit de montrer que  
$Lf^{*} : LSEns_{\mathbb{U}}/S \longrightarrow LSEns_{\mathbb{U}}/S'$ est  
une \'equivalence de $1$-cat\'egories de Segal. Pour cela il suffit de remarquer 
que l'adjonction $(f_{!},f^{*})$ est une \'equivalence de Quillen, ce qui est  
varie car $SEns$ est une cmf propre \`a droite (i.e. les changements de bases le long
de fibrations pr\'eservent les \'equivalences). \hfill $\Box$ \\ 
 
Le lemme implique en particulier que le foncteur  
$(Fib(-),L\omega)$ se factorise par la cat\'egorie  
homotopique,  
$$(Fib(-),L\omega) : Ho(*/SEns_{\mathbb{U}})\longrightarrow Ho(*/1-PrCat_{\mathbb{V}}).$$ 
 
\begin{prop}\label{p2} 
Pour tout $X \in Es^{cw}_{\mathbb{U}}$, il existe une adjonction de Quillen  
$$R : SPr(X) \longrightarrow SEns_{\mathbb{U}}/S(X) \qquad  
Se : SEns_{\mathbb{U}}/S(X) \longrightarrow SPr(X).$$ 
De plus, pour tout morphisme dans $Es_{\mathbb{U}}^{cw}$,  
$f : Y \longrightarrow X$, il existe des isomorphismes, naturels en $f$ 
$$f_{*} Se \simeq Se f_{*} \qquad R f^{*} \simeq f^{*} R.$$ 
\end{prop} 
 
\textit{Preuve:} Soit $X \in Es^{cw}_{\mathbb{U}}$ et d\'efinissons l'adjonction de Quillen de la fa\c{c}on suivante. 
Pour $Y \rightarrow S(X)$ un objet de $SEns_{\mathbb{U}}/S(X)$ on d\'efinit $Se(Y\rightarrow S(X))$ par 
$$\begin{array}{ccc} 
Ouv(X)^{o} & \longrightarrow & SEns_{\mathbb{U}} \\ 
(U\subset X) & \mapsto & \underline{Hom}_{S(X)}(S(U),Y) 
\end{array}$$ 
o\`u $\underline{Hom}_{S(X)}$ est le $Hom$ simplicial de $SEns_{\mathbb{U}}/S(X)$.  
 
Son adjoint \`a gauche est d\'efini de la fa\c{c}on suivante. Pour $F \in SPr(X)$ on forme le foncteur 
$$\begin{array}{ccc} 
Ouv(X)^{o}\times Ouv(X) & \longrightarrow & SEns_{\mathbb{U}}/S(X) \\ 
(U,V) & \mapsto & F(U)\times S(V) 
\end{array}$$ 
L'objet $R(F)$ est alors le co-end de ce foncteur dans $SEns_{\mathbb{U}}/S(X)$. Plus pr\'ecis\'ement, on a  
$$R(F):=\left( \coprod_{U \in Ouv(X)}F(U)\times S(U) \right)/\mathcal{R}$$ 
o\`u $\mathcal{R}$ identifie $(u^{*}(x),s)\in F(U)\times S(U)$ avec $(x,u_{*}(s))\in F(V)\times S(V)$, pour une inclusion 
$u : U \subset V$ et $(x,s) \in F(V)\times S(U)$. Il est alors facile de voir que $Se$ est adjoint \`a droite de $R$. \\ 
 
Commen\c{c}ons par remarquer que $Se$ est de Quillen \`a droite lorsque  
$Ouv(X)$ est muni de la topologie triviale. En d'autres termes il faut montrer  
que pour toute fibration (resp. fibration 
triviale) $Y \longrightarrow Y'$ dans $SEns_{\mathbb{U}}/S(X)$, et tout $U \in Ouv(X)$, le morphisme 
$$Se(Y)(U) \longrightarrow Se(Y')(U)$$ 
est une fibration (resp. une fibration triviale) dans $SEns_{\mathbb{U}}$.  
Mais ceci se d\'eduit imm\'ediatement du fait que $SEns_{\mathbb{U}}/S(X)$ est une cmf 
simpliciale et que tous ses objets sont cofibrants. 
  
Ainsi, pour montrer que $(R,Se)$ reste une adjonction de Quillen lorsque  
l'on effectue la localisation de Bousfield \`a gauche, il 
suffit de montrer que $Se$ pr\'eserve les objets fibrants. Or, d'apr\`es  
le lemme \ref{l4'} les objets fibrants dans $SPr(X)$ sont 
les pr\'efaisceaux $F$ fibrants pour la topologie triviale et qui  
satisfont de plus \`a la propri\'et\'e de descente 
pour les hyper-recouvrements. Il nous faut donc montrer que pour  
toute fibration $Y \rightarrow S(X)$, tout ouvert $U \subset X$, et tout 
hyper-recouvrement ouvert $U_{*}$ de $U$, le morphisme naturel 
$$\underline{Hom}_{S(X)}(S(U),Y) \longrightarrow Holim_{[m] \in \Delta}\underline{Hom}_{S(X)}(S(U_{m}),Y)$$
$$\simeq \underline{Hom}_{S(X)}(Hocolim_{[m] \in \Delta^{o}}S(U_{m}),Y)$$ 
est une \'equivalence. 
Or, ce morphisme est induit par le morphisme naturel 
$$u : Hocolim_{[m]\in \Delta^{o}}S(U_{m}) \longrightarrow S(Hocolim_{[m] \in \Delta^{o}}U_{m}) \longrightarrow S(U).$$ 
 
Rappelons alors le lemme suivant. 
 
\begin{lem} 
Le morphisme naturel  
$$Hocolim_{[m] \in \Delta^{o}}U_{m} \rightarrow U$$  
est une \'equivalence dans $Es_{\mathbb{U}}$. 
\end{lem} 
 
\textit{Preuve:} D'apr\'es les th\'eor\`emes $1.3.2$ et $1.3.5$ de \cite{fp} (voir aussi  
l'exercice $1$ du $\S 1.3$), tout ouvert d'un $CW$ complexe est paracompact et localement contractile. Le lemme 
est alors d\'emontr\'e au cours de la preuve de \cite[Thm. $2.1$]{am}.  \hfill $\Box$ \\ 
 
En utilisant le lemme pr\'ec\'edent et le fait que le foncteur $S$ pr\'es\`erve les colimites homotopiques,  
on voit que le morphisme $u$ est  
une \'equivalence d'ensembles simpliciaux.  Ainsi, comme $Y\rightarrow S(X)$  
est fibrant dans $SEns_{\mathbb{U}}/S(X)$, le morphisme en question 
$$\underline{Hom}_{S(X)}(S(U),Y) \longrightarrow Holim_{[m] \in \Delta}\underline{Hom}_{S(X)}(S(U_{m}),Y)$$ 
est une \'equivalence. Ceci montre que $Se$ pr\'eserve les objets  
fibrants, et donc que les foncteurs $R$ et $Se$ forment une  
adjonction de Quillen. \\ 
 
Soit $f : Y \longrightarrow X$, et $Z\rightarrow S(Y)$ dans  
$SEns_{\mathbb{U}}/S(Y)$. Alors par d\'efinition on a pour tout ouvert $U \subset X$ 
$$f_{*}Se(Z)(U)\simeq \underline{Hom}_{S(Y)}(S(U)\times_{S(X)}S(Y),Z)\simeq 
\underline{Hom}_{S(Y)}(f^{*}S(U),Z)$$
$$\simeq \underline{Hom}_{S(X)}(S(U),f_{*}(Z))=Se f_{*}(Z).$$ 
Ces isomorphismes \'etant fonctoriels en $U$ et $Z$, ils d\'eterminent  
un isomorphisme $f_{*}Se\simeq Se f_{*}$. Le lecteur v\'erifiera sans 
peine que cet isomorphisme est naturel en $f$. L'isomorphisme  
$R f^{*} \simeq f^{*} R$ se d\'eduit alors du pr\'ec\'edent par adjonction. \hfill $\Box$ \\ 
 
\begin{lem}\label{l1'} 
Soit $f : Y \longrightarrow X$ un morphisme dans $Es^{cw}_{\mathbb{U}}$.  
Alors, pour tout $Z \in Ho(SEns_{\mathbb{U}}/S(X))$ il existe un 
isomorphisme dans $Ho(SPr(Y))$, naturel en $Z$ 
$$\mathbb{R}Se\mathbb{R}f^{*}(Z)\simeq f^{*}\mathbb{R}Se(Z).$$ 
\end{lem} 
 
\textit{Preuve:} Par adjonction entre $f^{*}$ et $f_{*}$,  
l'isomorphisme $f_{*} Se \simeq Se f_{*}$ induit un morphisme naturel 
en $f$, $u : Se f^{*} \longrightarrow f^{*} Se$.  
 
Soit $Z \in Ho(SEns_{\mathbb{U}}/S(X))$. Alors  
$\mathbb{R}Se\mathbb{R}f^{*}(Z)\simeq Se f^{*}(Z')$, o\`u $Z \longrightarrow Z'$ est un mod\`ele 
fibrant dans $SEns_{\mathbb{U}}/S(X)$. En composant avec $u$, on trouve donc un morphisme naturel en $Z$ 
$$u_{Z} : \mathbb{R}Se\mathbb{R}f^{*}(Z) \longrightarrow f^{*}Se(Z')\simeq f^{*}\mathbb{R}Se(Z).$$ 
Soit $y \in Y$. La fibre en $y$ du morphisme ci-dessus est isomorphe dans $Ho(SEns_{\mathbb{U}})$ \`a  
$$Se(f^{*}(Z'))_{y} \longrightarrow (f^{*}Se(Z'))_{y}\simeq Se(Z')_{f(y)}.$$ 
Ainsi, pour montrer que $u_{Z}$ est un isomorphisme,  
il suffit de montrer que pour tout $p : S \rightarrow S(Y)$ fibrant dans $SEns_{\mathbb{U}}/S(Y)$,  
le morphisme naturel $Se(S)_{y} \longrightarrow p^{-1}(y)$  
est une \'equivalence. Par d\'efinition, on a  
$$Se(S)_{y}=Colim_{y\in V}\underline{Hom}_{S(Y)}(S(V),S) \longrightarrow \underline{Hom}_{S(Y)}(S(y),S)\simeq p^{-1}(y),$$ 
o\`u la colimite est prise sur les voisinages $V$ de $y$  
dans $Y$. Or, comme $Y$ est localement contractile et $p$ une fibration,  
le morphisme ci-dessus est une colimite filtrante  
d'\'equivalences dans $SEns_{\mathbb{U}}$, et est donc une \'equivalence. \hfill $\Box$ \\ 
 
\begin{lem}\label{l2'} 
Soit $f : Y \longrightarrow X$ une fibration dans $Es^{cw}_{\mathbb{U}}$.  
 Alors, pour tout $F \in Ho(SPr(X))$ il existe un 
isomorphisme dans $Ho(SEns_{\mathbb{U}}/S(Y))$, naturel en $F$ 
$$\mathbb{L}Rf^{*}(F)\simeq \mathbb{R}f^{*}\mathbb{L}R(F).$$ 
\end{lem} 
 
\textit{Preuve:} Consid\'erons l'adjoint \`a droite de $f^{*}$,  
$f_{*} : SEns_{\mathbb{U}}/S(Y) \longrightarrow SEns_{\mathbb{U}}/S(X)$. 
Comme $f$ est une fibration, le foncteur $f_{*}$ de Quillen \`a droite, 
 et donc $f^{*}$ est de Quillen \`a gauche. En particulier, comme 
$f^{*}$ pr\'eserve les \'equivalences, on a $\mathbb{R}f^{*}\simeq \mathbb{L}f^{*}$.  
Le lemme se d\'eduit alors de la derni\`ere assertion 
de la proposition \ref{p2}. \hfill $\Box$ \\ 
 
Nous pouvons maintenant \'enoncer et d\'emontrer le premier th\'eor\`eme de classification. 
 
\begin{thm}\label{p3} 
Pour tout $X \in Es^{cw}_{\mathbb{U}}$, le morphisme  
$$LSe : Fib(S(X)) \longrightarrow LSPr(X)$$ 
se factorise par la sous-$1$-cat\'egorie de Segal pleine de  
$LSPr(X)$ form\'ee des objets localement $h$-constants. De plus, le morphisme 
induit 
$$LSe : Fib(S(X)) \longrightarrow Loc(X)$$ 
est un un isomorphisme dans $Ho(1-PrCat_{\mathbb{V}})$, et son inverse est $LR$. 
\end{thm} 
 
\textit{Preuve:} Nous allons proc\'eder en plusieurs \'etapes. 
 
\begin{lem}\label{l2} 
Pour tout $X \in Es^{cw}_{\mathbb{U}}$, et tout $Y \in SEns_{\mathbb{U}}/S(X)$,  
l'objet $\mathbb{R}Se(Y) \in Ho(SPr(X))$ est localement $h$-constant. 
\end{lem} 
 
\textit{Preuve:} Comme l'assertion est locale sur $X$, on peut  
gr\^ace au lemme \ref{l1'} supposer que $X$ est contractile.  
L'ensemble simplicial $S(X)$  
se r\'etracte alors par d\'eformation sur un des ses points $s \in S(X)$. Soit  
$Y \longrightarrow S(X)$ fibrant dans $SEns_{\mathbb{U}}/S(X)$.  
Comme $Y \longrightarrow S(X)$ est une fibration, $Y \rightarrow S(X)$ est  
isomorphe dans $Ho(SEns_{\mathbb{U}}/S(X))$ \`a  
$Y_{s}\times S(X) \longrightarrow S(X)$, o\`u $Y_{s}$ est la fibre de $Y$ en $s$.  
Ainsi, $\mathbb{R}Se(Y)$ est isomorphe \`a  
$\mathbb{R}Se(Y_{s}\times S(X))\simeq Se(Y_{s}\times S(X))$.  
Mais $Se(Y_{s}\times S(X))$ est le pr\'efaisceau d\'efini par 
$$\begin{array}{cccc} 
Se(Y_{s}\times S(X)) : & Ouv(X)^{o} & \longrightarrow & SEns_{\mathbb{U}} \\ 
& (U \subset X) & \mapsto & \underline{Hom}(S(U),Y_{s}) 
\end{array}$$ 
Soit $\underline{Y_{s}}$ le pr\'efaisceau constant de  
fibre $Y_{s}$. On d\'efinit un morphisme dans $SPr(X)$,  
$u : \underline{Y_{s}} \longrightarrow Se(Y_{s}\times S(X))$,  
qui au-dessus de $U \subset X$ est le morphisme 
$Y_{s} \longrightarrow \underline{Hom}(S(U),Y_{s})$  
correspondant par adjonction \`a la projection $Y_{s}\times S(U) \longrightarrow Y_{s}$.  
Par le m\^eme argument que celui utilis\'e dans la preuve  
du lemme \ref{l1'}, on montre que les fibres de $u$ sont 
des colimites filtrantes d'\'equivalences, et donc des  
\'equivalences. Ceci implique donc que $Se(Y_{s}\times S(X))$ est  
$h$-constant et donc localement $h$-constant. \hfill $\Box$ \\ 
 
\textit{Remarque:} Lors de la preuve du lemme ci-dessus on a  
aussi montr\'e que pour tout $X \in Es^{cw}_{\mathbb{U}}$, et  
tout $Y \in SEns_{\mathbb{U}}$, le morphisme naturel  
$\underline{Y} \longrightarrow \mathbb{R}Se(Y\times S(X))$ 
\'etait un isomorphisme dans $Ho(SPr(X))$ (o\`u  
$\underline{Y}$ est le pr\'efaisceau constant de fibre $Y$). 
 
\begin{lem}\label{l3} 
Pour tout $X\in Es_{\mathbb{U}}^{cw}$ et tout  
$Y \in Ho(SEns_{\mathbb{U}}/S(X))$, le morphisme d'adjonction 
$$\mathbb{L}R\mathbb{R}Se(Y) \longrightarrow Y$$ 
est un isomorphisme.   
\end{lem} 
 
\textit{Preuve:}  Soit $x \in X$ un point de $X$, et  
$f : (P,p) \longrightarrow (X,x)$ une fibration, avec $P$ contractile.  
En appliquant $\mathbb{R}f^{*}$ au morphisme en  
question et un utilisant le lemme \ref{l2'} on obtient un morphisme 
dans $Ho(SEns_{\mathbb{U}}/S(P))$ 
$$\mathbb{R}f^{*}\mathbb{L}R\mathbb{R}Se(Y)\simeq \mathbb{L}R\mathbb{R}Se f^{*}(Y) \longrightarrow \mathbb{R}f^{*}(Y).$$ 
Comme $S(P)$ est contractile, $\mathbb{R}f^{*}(Y)\simeq Z\times S(P)$.  
On a alors d\'ej\`a vu lors de la preuve du lemme \ref{l2} que  
$\mathbb{R}Se\mathbb{R}f^{*}(Y)$ est naturellement  
isomorphe au pr\'efaisceau constant de fibre $Z$. Notons 
$\underline{Z}$ ce pr\'efaisceau. 
Or, comme $\underline{Z}$ est un objet cofibrant,  
$\mathbb{L}R\mathbb{R}Se\mathbb{R}f^{*}(Y)$ 
est isomorphe \`a $R(\underline{Z})\simeq Z\times S(P)$. A travers ces identifications, le morphisme  
$$\mathbb{R}f^{*}\mathbb{L}R\mathbb{R}Se(Y) \longrightarrow \mathbb{R}f^{*}(Y)$$ 
est isomorphe \`a l'identit\'e. C'est donc un isomorphisme. \\ 
 
On vient donc de montrer que pour tout $x \in X$,  
le morphisme induit par la fibration $f : (P,p) \longrightarrow (X,x)$ 
$$u : \mathbb{R}f^{*}\mathbb{L}R\mathbb{R}Se(Y) \longrightarrow \mathbb{R}f^{*}(Y)$$ 
est un isomorphisme dans $Ho(SEns_{\mathbb{U}}/S(P))$.  
Or, pour $Y \in Ho(SEns_{\mathbb{U}}/S(X))$, $\mathbb{R}f^{*}(Y)$ est isomorphe \`a la fibre  
homotopique de $Y$ en $x \in S(X)$.  
Ceci signifie que le morphisme  
d'adjonction  
$$\mathbb{L}R\mathbb{R}Se(Y) \longrightarrow Y$$ 
est un morphisme induisant une \'equivalence sur toutes  
les fibres homotopiques. C'est donc une \'equivalence dans  
$Ho(SEns_{\mathbb{U}}/S(X))$. \hfill $\Box$ \\ 
 
Remarquons d\`es \`a pr\'esent que les lemmes \ref{l2} et \ref{l3}  
impliquent que $\mathbb{R}Se : Ho(SEns_{\mathbb{U}}/S(X)) \longrightarrow Ho(Loc(X))$ 
est pleinement fid\`ele. Il en est donc de m\^eme du  
foncteur $\mathbb{R}Se : Ho(SEns_{\mathbb{U}}/S(X)) \longrightarrow Ho(SPr(X))$. De plus,  
le foncteur $Se$ \'etant compatible avec l'exponentiation  
(i.e. il existe des isomorphismes naturels 
$Se(Y^{A})\simeq Se(Y)^{A}$, pour tout  
$Y \in SEns_{\mathbb{U}}/S(X)$ et $A \in SEns_{\mathbb{U}}$), on en d\'eduit que les morphismes naturels dans $Ho(SEns_{\mathbb{U}})$ 
$$\mathbb{R}\underline{Hom}_{S(X)}(Y,Z) \longrightarrow \mathbb{R}\underline{Hom}_{X}(\mathbb{R}Se(Y),\mathbb{R}Se(Z)),$$ 
sont des isomorphismes dans $Ho(SPr(X))$, et ce pour  
tout $Y,Z \in Ho(SEns_{\mathbb{U}}/S(X))$. Cette derni\`ere remarque nous sera utile pour d\'emontrer 
le lemme suivant.   
 
\begin{lem}\label{l4} 
Pour tout $X \in Es_{\mathbb{U}}^{cw}$ et tout objet $F \in Ho(Loc(X))$ le morphisme d'adjonction  
$F \longrightarrow \mathbb{R}Se\mathbb{L}RF$ est un isomorphisme. 
\end{lem} 
 
\textit{Preuve:} Choisissons une fibration surjective $f : P \longrightarrow X$, avec $P$ contractile.  
Le corollaire \ref{c3'} et les lemmes \ref{l1'} et \ref{l2'}  
montrent que l'on peut alors supposer que $X=P$. On supposera donc que $X$ est contractile. \\ 
 
Commen\c{c}ons par montrer que tout $F \in PrLoc(X)$ est  
alors $h$-constant. Nous pouvons sans perte de g\'en\'eralit\'e supposer 
que $F$ est fibrant. En choisissant une r\'etraction de $X$  
sur l'un de ses points on v\'erifie qu'il suffit de montrer que pour tout $F \in PrLoc(X\times I)$, le morphisme d'adjonction 
$$\mathbb{L}p^{*}\mathbb{R}p_{*}(F) \longrightarrow F$$ 
est un isomorphisme, ou $p : X\times I \longrightarrow X$  
est la premi\`ere projection, et  
$I=[0,1]$ est l'intervalle standard. Il est alors facile  
de voir qu'il est \'equivalent de montrer pour tout $F$ fibrant dans $PrLoc(I)$, $F$ est $h$-constant,  
et que et le morphisme fibre en $0 \in I$, $F(I) \longrightarrow F_{0}$,  
est une \'equivalence dans $SEns_{\mathbb{U}}$.  \\ 
 
Soit $t$ la borne sup\'erieure des $x \in I$ tel que la  
restriction de $F$ \`a $[0,x[$ soit $h$-constant. Supposons que $t<1$ et soit 
$x>0$ tel que la restriction de $F$ \`a $]t-x,t+x[$ soit  
$h$-constant. Notons alors $Z_{0}$ et $Z_{1}$ deux ensembles simpliciaux fibrants 
dans $SEns_{\mathbb{U}}$, et $u_{0} : \underline{Z_{0}}\simeq F_{|[0,t[}$,  
$u_{1} : \underline{Z_{1}} \simeq F_{|]t-x,t+x[}$ des 
isomorphismes dans la cat\'egorie homotopique (o\`u  
$\underline{Z}$ est le pr\'efaisceau constant de fibre $Z$).  
Remarquons que les objets $\underline{Z_{0}}$ et $\underline{Z_{1}}$ sont  
cofibrants. Comme $F$ est fibrant, ses restrictions  
aussi, et on peut ainsi repr\'esenter les morphismes pr\'ec\'edents par  
des morphismes de pr\'efaisceaux simpliciaux 
$$u_{0} : \underline{Z_{0}} \longrightarrow F_{|[0,t[}, \qquad  
u_{1} : \underline{Z_{1}} \longrightarrow F_{|]t-x,t+x[}.$$ 
De plus, l'isomorphisme $\beta:=u_{1}^{-1}u_{0} : \underline{Z_{0}} \longrightarrow \underline{Z_{1}}$,  
est un isomorphisme entre 
$\mathbb{R}Se(Z_{0}\times S(]t-x,t[))$ et $\mathbb{R}Se(Z_{1}\times S(]t-x,t[))$.  
Or, d'apr\`es le lemme \ref{l3} le foncteur 
$\mathbb{R}Se$ est pleinement fid\`ele. Ainsi, comme les objets  
$Z_{0}\times S(]t-x,t[)$ et $Z_{1}\times S(]t-x,t[)$ sont 
cofibrants et fibrants dans $SEns_{\mathbb{U}}/S(]t-x,t[)$,  
on peut trouver un diagramme commutatif \`a homotopie pr\`es dans $SPr(]t-x,t[)$ 
$$\xymatrix{ 
\underline{Z_{0}} \ar[dr]_-{u_{0}} \ar[rr]^-{\beta} & & \underline{Z_{1}} \ar[dl]^-{u_{1}} \\ 
& F_{|]t-x,t[}. & }$$ 
De plus, quitte \`a bouger $u_{0}$ dans sa classe d'homotopie  
on peut supposer que ce diagramme est strictement commutatif. On d\'efinit alors 
$$u : \underline{Z_{0}}  \longrightarrow  F_{|[0,t+x[}$$ 
comme \'etant $u_{0}$ sur $[0,t[$ et $u_{1}\beta$ sur $]t-x,t+x[$.  
Comme $u_{0}$ et $u_{1}$ sont des \'equivalences $u$ est une \'equivalence 
dans $SPr([0,t+x[)$. Ceci contredit $t<1$, et donc $t=1$. 
 L'objet $F \in SPr(I)$ est donc $h$-constant.    
 
On peut alors \'ecrire $F \simeq \mathbb{R}Se(Z\times S(I))$,  
o\`u $Z \in SEns_{\mathbb{U}}$. Or, comme on sait que $\mathbb{R}Se$ est pleinement fid\`ele,  
le morphisme fibre en $0$, $F(I)\simeq \mathbb{R}\underline{Hom}(*,F) \longrightarrow F_{0}$  
est isomorphe dans $Ho(SEns_{\mathbb{U}})$ au morphisme 
$$\mathbb{R}\underline{Hom}_{S(I)}(S(I),Z\times S(I)) \longrightarrow Z\times \{0\},$$ 
qui \`a $f : S(I) \longrightarrow Z \times S(I)$ associe  
$f(0)$. Mais ce dernier morphisme est clairement un isomorphisme. Ceci termine 
donc la preuve que tout $F \in SPr(X)$ est $h$-constant. \\ 
 
Revenons donc au cas o\`u $X$ est contractile, et $F \in PrLoc(X)$. 
 Comme on sait que $F$ est $h$-constant on peut supposer que $F$ est constant 
de fibre $Z \in SEns_{\mathbb{U}}$, que l'on supposera fibrant.  
Dans ce cas, $F$ est cofibrant et $\mathbb{R}L(F)\simeq R(F)\simeq Z\times S(X)$. De plus, comme 
$Z$ est fibrant, $Z\times S(X)$ est fibrant dans  
$SEns_{\mathbb{U}}/S(X)$. Ainsi, $\mathbb{R}Se(F)\simeq \underline{Z}$, et  
le morphisme $\mathbb{L}R\mathbb{R}Se(F) \longrightarrow F$  
est isomorphe au morphisme naturel 
$$R(\underline{Z})=R(Z\times *) \longrightarrow Z\times R(*)\simeq Z\times S(X),$$ 
qui est un isomorphisme. \hfill $\Box$ \\ 
 
Soit $C=PrLoc(X)^{cf}$ la sous-cat\'egorie pleine de $SPr(X)$  
form\'ee des objets fibrants, cofibrants et localement $h$-constants. De m\^eme, soit 
$D=(SEns_{\mathbb{U}}/S(X))^{cf}$ la sous-cat\'egorie pleine  
de $SEns_{\mathbb{U}}/S(X)$ form\'ee des objets fibrants et cofibrants. Nous noterons 
$Q \rightarrow Id$ (resp. $Id \rightarrow P$) les foncteurs de remplacement cofibrants (resp. remplacement fibrants) dans 
$SPr(X)$ ou $SEns_{\mathbb{U}}/S(X)$ (voir \cite[$\S 1.1$]{ho}). 
 On d\'efinit deux foncteurs $R'$ et $Se'$ par 
$$R':=P\circ R : C \longrightarrow D \qquad Se':=Q\circ Se : D \longrightarrow C.$$ 
On dispose alors de diagrammes de foncteurs 
$$\xymatrix{R'Se'=PRQSe \ar[r]^-{a} & P & \ar[l]_-{b} Id}$$ 
$$\xymatrix{Se'R'=QSePR & \ar[l]_-{c}  Q  \ar[r]^-{d} & Id}.$$ 
De plus, les transformations naturelles $b$ et $d$ sont 
 des \'equivalences objet par objet. De m\^eme, les lemmes \ref{l3} et \ref{l4}  
impliquent que $a$ et $c$ sont des \'equivalences objets par objets. Ceci implique que les morphismes $LR' : LC \longrightarrow LD$  
et $LSe' : LD \longrightarrow LC$ sont inverses l'un de  
l'autre dans $Ho(1-PrCat_{\mathbb{V}})$, et ach\`eve donc la preuve  
du th\'eor\`eme \ref{p3}. \hfill $\Box$ \\ 
 
\begin{cor}\label{c2} 
Pour tout $(X,x) \in */Es^{cw}_{\mathbb{U}}$ il existe un  
isomorphisme dans $Ho(*/1-PrCat_{\mathbb{V}})$, fonctoriel en $(X,x)$ 
$$(\mathbb{R}\underline{Hom}(Loc(X),Top),L\omega_{x}) \simeq (\mathbb{R}\underline{Hom}(Fib(S(X)),Top),L\omega_{x}).$$ 
\end{cor} 
 
\textit{Preuve:} Soit $p : (P,y) \longrightarrow (X,x)$ une  
fibration dans $Es_{\mathbb{U}}^{cw}$, avec $P$ contractile. 
Par la proposition \ref{p2} il existe un carr\'e de cat\'egories  
$$\xymatrix{ 
PrLoc(X)^{c} \ar[r]^-{R} \ar[d]_-{p^{*}} & SEns_{\mathbb{U}}/S(X) \ar[d]^-{p^{*}} & \\ 
PrLoc(P)^{c} \ar[r]_-{R} & SEns_{\mathbb{U}}/S(P)}$$  
et un isomorphisme naturel de foncteurs $h : p^{*}R\simeq Rp^{*}$.  
On repr\'esente cet isomorphisme par un diagramme commutatif 
de cat\'egories 
$$\xymatrix{PrLoc(X)^{c} \ar[rr]^-{R} \ar[d]_-{i_{0}} & & SEns_{\mathbb{U}}/S(X) \ar[d]^-{p^{*}} \\ 
PrLoc(X)^{c}\times \overline{I} \ar[rr]^-{h} & & SEns_{\mathbb{U}}/S(P) \\ 
PrLoc(X)^{c} \ar[u]^-{i_{1}} \ar[rr]_-{p^{*}} & & PrLoc(P)^{c}, \ar[u]_-{R} }$$ 
o\`u $\overline{I}$ est la cat\'egorie poss\'edant  
deux objets $0$ et $1$ et un unique isomorphisme entre les deux.  
En appliquant la localisation de Dwyer-Kan on obtient  
un diagramme commutatif dans $1-PrCat_{\mathbb{V}}$ 
$$\xymatrix{LPrLoc(X)^{c} \ar[rr]^-{LR} \ar[d]_-{Li_{0}} & & Fib(S(X)) \ar[d]^-{Lp^{*}} \\ 
L(PrLoc(X)^{c}\times \overline{I}) \ar[rr]^-{Lh} & & Fib(S(P)) \\ 
LPrLoc(X)^{c} \ar[u]^-{Li_{1}} \ar[rr]_-{Lp^{*}} & & LPrLoc(P)^{c}. \ar[u]_-{LR} }$$ 
Par le th\'eor\`eme \ref{p3} les morphismes $LR$  
sont des \'equivalences. De plus, comme $i_{0}$ et $i_{1}$ sont des \'equivalences de 
cat\'egories, $Li_{0}$ et $Li_{1}$ sont aussi des \'equivalences.  
Ainsi, ce diagramme fournit un isomorphisme bien 
d\'efini dans $Ho(1-PrCat)$ 
$$(\mathbb{R}\underline{Hom}(LPrLoc(X)^{c},LPrLoc(P)^{c}),Lp^{*})) 
\longrightarrow (\mathbb{R}\underline{Hom}(Fib(S(X)),Fib(S(P))),Lp^{*})).$$ 
De plus, une application du lemme \ref{l1}, de  
la proposition \ref{p2} et du th\'eor\`eme \ref{p3} implique qu'il existe des isomorphismes naturels 
$$(\mathbb{R}\underline{Hom}(LPrLoc(X)^{c},LPrLoc(P)^{c}),Lp^{*}))\simeq  
(\mathbb{R}\underline{Hom}(Loc(X),Top),L\omega_{x})$$ 
$$(\mathbb{R}\underline{Hom}(Fib(S(X)),Fib(S(P))),Lp^{*})) \simeq   
(\mathbb{R}\underline{Hom}(Fib(S(X)),Top),L\omega_{x})$$ 
On obtient ainsi l'isomorphisme cherch\'e.  
En utilisant que l'isomorphisme $h$ est naturel en $f$, Le lecteur v\'erifiera que cet 
isomorphisme est fonctoriel en $(X,x) \in */Es^{cw}_{\mathbb{U}}$. \hfill $\Box$ \\ 
 
Remarquons enfin que le corollaire \ref{c2} et le lemme \ref{l1}  
impliquent que le foncteur $(X,x) \mapsto  (\mathbb{R}\underline{Hom}(Loc(X),Top),L\omega_{x})$ induit un 
foncteur sur la cat\'egorie homotopique 
$$\begin{array}{ccc} 
Ho(*/Es_{\mathbb{U}}^{cw}) & \longrightarrow & Ho(*/1-PrCat_{\mathbb{V}}) \\ 
(X,x) & \mapsto & (\mathbb{R}\underline{Hom}(Loc(X),Top),L\omega_{x}). 
\end{array}$$ 
De plus, comme tout espace est (faiblement) \'equivalent \`a un  
$CW$-complexe, le foncteur naturel $Ho(*/Es^{cw}_{\mathbb{U}})  
\longrightarrow Ho(*/Es_{\mathbb{U}})$ est une 
\'equivalence de cat\'egories. Ainsi, le foncteur pr\'ec\'edent  
induit un foncteur bien d\'efini \`a isomorphisme unique pr\`es 
$$\begin{array}{ccc} 
Ho(*/Es_{\mathbb{U}}) & \longrightarrow & Ho(*/1-PrCat_{\mathbb{V}}) \\ 
(X,x) & \mapsto & (\mathbb{R}\underline{Hom}(Loc(X),Top),L\omega_{x}). 
\end{array}$$ 
Notons que ceci est un abus de notations car l'image de  
$(X,x) \in Ho(*/Es_{\mathbb{U}})$ est r\'eellement d\'efinie comme \'etant 
$(\mathbb{R}\underline{Hom}(Loc(X'),Top),L\omega_{x'})$,  
o\`u $(X',x')$ est un $CW$ complexe point\'e (faiblement) \'equivalent \`a $(X,x)$. 
  
\end{subsection} 
 
\begin{subsection}{Fibrations et repr\'esentations des groupes simpliciaux} 
 
Dans ce paragraphe nous allons d\'emontrer un analogue pour le  
cas des fibrations de la correspondance entre rev\^etements sur $X$ et  
$\pi_{1}(X,x)$-ensembles. Une remarque cl\'e est que tout  
ensemble simplicial $S$ connexe est \'equivalent au  
classifiant d'un groupe simplicial $G$. Ainsi, d'apr\`es la  
propri\'et\'e d'invariance de la $1$-cat\'egorie de Segal des fibrations 
sur $S$ (voir le lemme \ref{l1}), on supposera que $S=d(G)$  
est un tel classifiant. On d\'efinira alors une \'equivalence de Quillen entre 
$SEns/S$ et la cat\'egorie des $G$-ensembles simpliciaux de  
fa\c{c}on assez \'evidente. Ceci montrera que \textit{les fibrations sur $d(G)$  
sont classifi\'es par les repr\'esentations de $G$}. \\ 
 
Soit $G \in SGp_{\mathbb{U}}$ un groupe simplicial et 
$G-SEns_{\mathbb{U}}$ la cat\'egorie des $G$-ensembles simpliciaux de $\mathbb{U}$ (i.e. 
des ensembles simpliciaux de $\mathbb{U}$ 
munis d'une action de $G$).  
On dispose d'une paire de foncteurs adjoints 
$$F : SEns_{\mathbb{U}} \longrightarrow G-SEns_{\mathbb{U}} \qquad  
\omega_{G} : G-SEns_{\mathbb{U}} \longrightarrow SEns_{\mathbb{U}}$$ 
o\`u $\omega_{G}$ est le foncteur d'oubli de la structure de  
$G$-modules. Le foncteur $F$ (\textit{F} pour \textit{free})  
associe \`a tout $T \in SEns_{\mathbb{U}}$ le $G$-module 
libre engendr\'e par $T$, d\'efini par $L(T):=T\times G$, o\`u  
$G$ op\`ere trivialement sur le facteur de gauche et par 
translations (\`a gauche) sur celui de droite. Pour plus de  
simplicit\'e nous noterons encore par $T$ l'ensemble simplicial $\omega_{G}(T)$, sous-jacent  
\`a un objet $T \in SEns_{\mathbb{U}}$. \\ 
 
Il existe une structure de cmf sur $G-SEns_{\mathbb{U}}$  
o\`u les fibrations (resp. les \'equivalences) sont les morphismes $f : T \longrightarrow T'$ 
tels que le morphisme induit sur les ensembles simpliciaux sous-jacents soit une fibration (resp. 
une \'equivalence). On peut aussi d\'efinir une  
structure de cat\'egorie enrichie sur $SEns_{\mathbb{U}}$ qui est fait une cat\'egorie  
de mod\`eles simpliciale (voir \cite[$\S V.2.3$]{gj}).  
 
Pour fixer les notations, rappelons que le produit  
externe de $T \in G-SEns_{\mathbb{U}}$ par un ensemble simplicial 
$S \in SEns_{\mathbb{U}}$ est simplement d\'efini par 
$S\otimes T:=S\times T$, o\`u $G$ op\`ere trivialement sur le premier facteur. Le 
$Hom$ simplicial entre deux $G$-ensembles simpliciaux $T$ et $T'$ sera not\'e 
$\underline{Hom}_{G}(T,T')$ et satisfait \`a la formule d'adjonction usuelle 
$$Hom(S\otimes T,T')\simeq Hom(S,\underline{Hom}_{G}(T,T')).$$ 
 
\begin{cor}\label{c3} 
Le foncteur $\omega_{G} : G-SEns_{\mathbb{U}} \longrightarrow SEns_{\mathbb{U}}$  
est de Quillen \`a droite.  
Le foncteur $F : SEns_{\mathbb{U}} \longrightarrow G-SEns_{\mathbb{U}}$ est de Quillen \`a gauche. 
\end{cor} 
 
\textit{Preuve:} C'est clair. \hfill $\Box$ \\ 
 
\begin{cor}\label{c4} 
Le foncteur 
$$F : SEns_{\mathbb{U}} \longrightarrow G-SEns_{\mathbb{U}}$$ 
est de Quillen \`a gauche, pr\'eserve la structure simpliciale et les objets fibrants. 
\end{cor} 
 
\textit{Preuve:} En effet, l'ensemble simplicial sous-jacent \`a $G$, est fibrant  
(\cite[$\S I.3.4$]{gj}). Ainsi,   
si $S \in SEns_{\mathbb{U}}$ est fibrant $\omega_{G}F(S)=S\times G$ est fibrant dans $SEns_{\mathbb{U}}$.  
Par d\'efinition ceci implique 
que $F(S)$ est fibrant dans $G-SEns_{\mathbb{U}}$. Il est clair  
que $F$ pr\'eserve les produits externes par des ensembles simpliciaux. \hfill $\Box$ \\ 
 
\begin{df}\label{d7'} 
Soit $G \in SGp_{\mathbb{U}}$. La $1$-cat\'egorie de Segal des repr\'esentations de $G$, $G-Top$, est d\'efinie par 
$$G-Top:=L(G-SEns_{\mathbb{U}}) \in 1-PrCat_{\mathbb{V}}.$$  
\end{df} 
 
Le foncteur $\omega_{G}$ pr\'eserve les \'equivalences et d\'efinit  
donc un morphisme dans $1-PrCat_{\mathbb{V}}$ 
$$L\omega_{G} :G-Top \longrightarrow Top.$$ 
 
Soit $f : G \longrightarrow H$ un morphisme dans $SGp_{\mathbb{U}}$. On en d\'eduit un foncteur 
$$f^{*} : H-SEns_{\mathbb{U}} \longrightarrow G-SEns_{\mathbb{U}}$$ 
obtenu en composant les actions par $f$. Ce foncteur  
pr\'eserve \'evidemment les \'equivalences, et de plus commute avec les foncteurs d'oublis. On 
dispose donc d'un diagramme commutatif  
$$\xymatrix{ 
H-SEns_{\mathbb{U}} \ar[rr]^-{f^{*}} \ar[dr]_-{\omega_{H}} & & G-SEns_{\mathbb{U}} \ar[dl]^-{\omega_{G}} \\ 
& SEns_{\mathbb{U}} & }$$ 
On y applique alors la localisation de Dwyer-Kan pour obtenir un morphisme dans $Ho(*/1-PrCat_{\mathbb{V}})$ 
$$(\mathbb{R}\underline{Hom}(G-Top,Top),L\omega_{G}) \longrightarrow (\mathbb{R}\underline{Hom}(H-Top,Top),L\omega_{H}).$$ 
Ceci d\'efini donc un foncteur 
$$\begin{array}{ccc} 
SGp_{\mathbb{U}} & \longrightarrow & Ho(*/1-PrCat_{\mathbb{V}}) \\ 
G & \mapsto & (\mathbb{R}\underline{Hom}(G-Top,Top),L\omega_{G}) 
\end{array}$$ 
 
\begin{lem} 
Soit $f : G \longrightarrow H$ une \'equivalence dans  
$SGp_{\mathbb{U}}$. Alors le morphisme induit dans $Ho(*/1-PrCat_{\mathbb{V}})$ 
$$(\mathbb{R}\underline{Hom}(G-Top,Top),L\omega_{G})  
\longrightarrow (\mathbb{R}\underline{Hom}(H-Top,Top),L\omega_{H})$$ 
est un isomorphisme. 
\end{lem} 
 
\textit{Preuve:} Il suffit de montrer que le foncteur  
$f^{*} : H-SEns_{\mathbb{U}} \longrightarrow G-SEns_{\mathbb{U}}$ est une \'equivalence de Quillen.  
Pour cela remarquons tout d'abord qu'il poss\`ede un adjoint \`a gauche 
$$\begin{array}{cccc} 
f_{!} : & G-SEns_{\mathbb{U}} & \longrightarrow & H-SEns_{\mathbb{U}} \\ 
& T & \mapsto & (T\times H)/G 
\end{array}$$ 
o\`u $G$ op\`ere diagonalement sur $T\times H$, et  
$H$ op\`ere sur $(T\times H)/G$ en op\'erant \`a droite sur le second facteur.  
Comme le foncteur $f^{*}$ est clairement de Quillen  
\`a droite, $f_{!}$ est de Quillen \`a gauche.  
 
Remarquons aussi que le foncteur  
$\mathbb{R}f^{*} : Ho(H-SEns_{\mathbb{U}}) \longrightarrow Ho(G-SEns_{\mathbb{U}})$  
est conservatif. Ainsi, il nous suffit  
de montrer que le morphisme d'adjonction  
$T \longrightarrow \mathbb{R}f^{*}\mathbb{L}f_{!}(T)$ est un isomorphisme 
pour tout $T \in Ho(G-SEns_{\mathbb{U}})$.  
Mais ce morphisme est isomorphe au morphisme naturel 
$$(T'\times G)/G \longrightarrow (T'\times H)/G,$$ 
o\`u $T'$ est un mod\`ele cofibrant pour $T$. Or,  
une application de \cite[$\S V.2.10$]{gj} montre que $T'\times G$ et $T'\times H$ sont 
cofibrants dans $G-SEns_{\mathbb{U}}$. Ainsi, comme  
le foncteur $X \mapsto X/G$ est de Quillen \`a gauche, et que le morphisme 
$T'\times G \longrightarrow T'\times H$ est une  
\'equivalence par hypoth\`ese, on en d\'eduit que  
$(T'\times G)/G \longrightarrow 
(T'\times H)/G$ est une \'equivalence d'ensembles  
simpliciaux. \hfill $\Box$ \\ 
 
Le lemme pr\'ec\'edent permet de d\'efinir un foncteur sur les cat\'egories homotopiques 
$$\begin{array}{ccc} 
H(SGp_{\mathbb{U}}) & \longrightarrow & Ho(*/1-PrCat_{\mathbb{V}}) \\ 
G & \mapsto & (\mathbb{R}\underline{Hom}(G-Top,Top),L\omega_{G}) 
\end{array}$$ 
 
Rappelons que pour $G \in SGp_{\mathbb{U}}$ on dispose  
de sa diagonale $d(G)$, qui est l'ensemble simplicial d\'efinit par  
$d(G)([m]):=G([m])^{m}$. De plus, comme $d(G)([0])=*$,  
cet ensemble simplicial est r\'eduit et donc naturellement point\'e. 
Ceci d\'efinit donc un foncteur  
$d : Ho(SGp_{\mathbb{U}}) \longrightarrow Ho(*/SEns_{\mathbb{U}})$. 
En le composant avec le foncteur  
$(S,s) \mapsto (\mathbb{R}\underline{Hom}(Fib(S),Top),L\omega_{s})$  
d\'efini dans le paragraphe pr\'ec\'edent on  en d\'eduit 
$$\begin{array}{ccc} 
Ho(SGp_{\mathbb{U}}) & \longrightarrow & Ho(*/1-PrCat_{\mathbb{V}}) \\ 
G & \mapsto & (\mathbb{R}\underline{Hom}(Fib(d(G)),Top),L\omega_{*}), 
\end{array}$$ 
o\'u $* \in d(G)$ d\'esigne le point de base naturel de la diagonale $d(G)$. 
 
Le second th\'eor\`eme de classification est le suivant. 
 
\begin{thm}\label{p5} 
Pour tout $G \in SGp_{\mathbb{U}}$ il existe un isomorphisme  
dans $Ho(*/1-PrCat_{\mathbb{V}})$ 
$$(\mathbb{R}\underline{Hom}(Fib(d(G)),Top),L\omega_{*}) \simeq  
(\mathbb{R}\underline{Hom}(G-Top,Top),L\omega_{G}),$$ 
qui est fonctoriel en $G \in Ho(SGp_{\mathbb{U}})$. 
\end{thm} 
 
\textit{Preuve:} Nous allons d\'efinir l'isomorphisme en question,  
et nous laisserons le soin au lecteur de v\'erifier que la construction 
est fonctorielle en $G$. \\ 
 
Rappelons que pour $G \in SGp_{\mathbb{U}}$ nous avons not\'e  
$EG$ la diagonale de l'ensemble bi-simplicial d\'efini par 
$$\begin{array}{cccc} 
EG : & \Delta^{o} & \longrightarrow & SEns_{\mathbb{U}} \\ 
& [m] & \mapsto & G^{m+1} 
\end{array}$$ 
et o\`u les faces et d\'eg\'en\'erescences sont  
donn\'ees par les projections et les diagonales. Comme le morphisme 
$G \rightarrow *$ poss\`ede une section, $EG$ est  
un ensemble simplicial contractile. De plus, il est naturellement muni d'une action de 
$G$ \`a gauche, et le quotient $EG/G$ est naturellement  
isomorphe \`a $d(G)$. On dispose donc d'un morphisme naturel  
dans $SEns_{\mathbb{U}}$ 
$$f_{G} : EG \longrightarrow d(G),$$ 
qui d'apr\`es \cite[$\S V.2.7$]{gj} est une fibration surjective.  
 
Notons $e \in EG$ le point correspondant \`a l'identit\'e dans  
le groupe $EG([0])=G([0])$. L'image de $e$ est \'evidemment le  
point de base $* \in d(G)$. On dispose donc 
d'un diagramme commutatif de foncteurs de Quillen \`a droite,  
$$\xymatrix{ 
SEns_{\mathbb{U}}/d(G) \ar[rd]_-{\omega_{\bullet}} \ar[rr]^-{f_{G}^{*}} & & SEns_{\mathbb{U}}/EG \ar[dl]^-{\omega_{e}} \\ 
& SEns_{\mathbb{U}} & }$$ 
En passant aux localis\'ees de Dwyer-Kan on en  
d\'eduit un morphisme dans $Ho(*/1-PrCat_{\mathbb{V}})$,  
$$L\omega_{e} : (\mathbb{R}\underline{Hom}(Fib(d(G)),Fib(EG)),Lf_{G}^{*}) \longrightarrow  
(\mathbb{R}\underline{Hom}(Fib(d(G)),Top),L\omega_{*}).$$ 
Or, d'apr\`es le lemme \ref{l1} le morphisme  
$L\omega_{e} : Fib(EG) \longrightarrow Top$ est un isomorphisme dans $Ho(1-PrCat)$. Ceci implique  
que $L\omega_{e}$ induit en r\'ealit\'e un isomorphisme  
$$L\omega_{e} : (\mathbb{R}\underline{Hom}(Fib(d(G)),Fib(EG)),Lf_{G}^{*}) \simeq  
(\mathbb{R}\underline{Hom}(Fib(d(G)),Top),L\omega_{*}).$$ 
Il nous suffit donc de montrer que le membre de gauche  
est isomorphe \`a $(\mathbb{R}\underline{Hom}(G-Top,Top),L\omega_{G})$. \\ 
 
Soit $G \in SGp_{\mathbb{U}}$ et d\'efinissons une adjonction de  
Quillen ($M$ pour \textit{monodromie}, et $D$ pour \textit{descente}) 
$$D : G-SEns_{\mathbb{U}} \longrightarrow SEns_{\mathbb{U}}/d(G)  
\qquad M : SEns_{\mathbb{U}}/d(G) \longrightarrow G-SEns_{\mathbb{U}}.$$ 
Pour $Y \rightarrow d(G)$ un objet de $SEns_{\mathbb{U}}/d(G)$ on pose 
$$M(Y):=\underline{Hom}_{d(G)}(EG,Y),$$ 
o\`u $G$ op\`ere sur $M(Y)$ en op\'erant sur $EG$. Son adjoint \`a gauche est d\'efini par  
$$D(T):=T\times_{G}EG=(T\times EG)/G,$$ 
pour tout $T \in G-SEns_{\mathbb{U}}$, et pour l'action  
diagonale sur $T\times EG$. D'apr\`es les propri\'et\'es des $Hom$ simpliciaux de $SEns_{\mathbb{U}}/d(G)$, 
il est imm\'ediat que $M$ est de Quillen \`a droite, et donc que $D$ est de Quillen \`a gauche.  
 
\begin{lem}\label{l5} 
Les foncteurs $M$ et $D$ d\'efinissent une \'equivalence de Quillen.  
\end{lem} 
 
\textit{Preuve:} Soit $Y \in SEns_{\mathbb{U}}/d(G)$ un objet fibrant. Alors, l'ensemble simplicial sous-jacent 
\`a $M(Y)$ est $\underline{Hom}_{d(G)}(EG,Y)$.  
Mais comme $EG$ est contractile, le morphisme naturel  
$$\underline{Hom}_{d(G)}(EG,Y) \longrightarrow \underline{Hom}_{d(G)}(e,Y)=Y_{e},$$ 
o\`u $Y_{e}$ est la fibre de $Y \rightarrow d(G)$ au  
point $*$, est une \'equivalence. Comme $d(G)$ est connexe ceci implique que le 
foncteur $\mathbb{R}M : Ho(SEns_{\mathbb{U}}/d(G)) \longrightarrow  
Ho(G-SEns_{\mathbb{U}})$ est conservatif. Il nous suffit donc de montrer que pour tout 
$T \in G-SEns_{\mathbb{U}}$ cofibrant le morphisme d'adjonction 
$$T \longrightarrow \mathbb{R}M\mathbb{L}D(T)$$ 
est un isomorphisme dans $Ho(G-SEns_{\mathbb{U}})$.  
Mais par d\'efinition de $M$ et $D$ ce morphisme est isomorphe \`a 
$$T \longrightarrow \mathbb{R}\underline{Hom}_{d(G)}(EG,T\times_{G}EG) 
\simeq \mathbb{R}\underline{Hom}_{EG}(EG,(T\times_{G}EG)\times_{d(G)}EG).$$ 
Or, $(T\times_{G}EG)\times_{d(G)}EG$ est naturellement  
isomorphe \`a $T\times EG$, et donc le morphisme pr\'ec\'edent se r\'eduit au morphisme 
naturel 
$$T \longrightarrow \mathbb{R}\underline{Hom}_{EG}(EG,T\times EG)\simeq \mathbb{R}\underline{Hom}(EG,T),$$ 
qui est clairement un isomorphisme car $EG$ est contractile. \hfill $\Box$ \\ 
 
Consid\'erons le diagramme suivant de cmf 
$$\xymatrix{ 
G-SEns_{\mathbb{U}} \ar[r]^-{D} \ar[d]_-{\omega_{G}} & SEns_{\mathbb{U}}/d(G) \ar[d]^-{f_{G}^{*}} \\ 
SEns_{\mathbb{U}} \ar[r]_-{-\times EG} & SEns_{\mathbb{U}}/EG}$$ 
o\`u $f_{G} : EG \longrightarrow d(G)$ est la projection naturelle.  
Comme ce diagramme commute \`a isomorphisme naturel pr\`es, 
et que les morphismes $LD$ et $L(-\times EG)$ sont des \'equivalences,  
le m\^eme argument que dans la preuve du corollaire \ref{c2} 
implique qu'il existe un isomorphisme fonctoriel en $G$ 
$$(\mathbb{R}\underline{Hom}(G-Top,Top),L\omega_{G}) \longrightarrow   
(\mathbb{R}\underline{Hom}(Fib(d(G)),Fib(EG)),Lf_{G}^{*}).$$ 
Enfin, nous avons d\'ej\`a vu  que le membre de droite est  
naturellement isomorphe \`a $(\mathbb{R}\underline{Hom}(Fib(d(G)),Top),L\omega_{*})$.  
\hfill $\Box$ \\ 
 
\end{subsection} 
 
\end{section} 
 
\begin{section}{Th\'eor\`eme de reconstruction} 
 
Soit $Ho(*/Es_{\mathbb{U}})_{c}$ la sous-cat\'egorie pleine de  
$Ho(*/Es_{\mathbb{U}})$ des espaces connexes par arcs. Nous disposons du foncteur 
(d\'efini \`a la fin du $\S 2.2$) 
$$\begin{array}{ccc} 
Ho(*/Es_{\mathbb{U}})_{c} & \longrightarrow & Ho(*/1-PrCat_{\mathbb{V}}) \\ 
(X,x) & \mapsto & (\mathbb{R}\underline{Hom}(Loc(X),Top),L\omega_{x}). 
\end{array}$$ 
On le compose avec $\mathbb{R}\Omega_{*} : Ho(*/1-PrCat_{\mathbb{V}})  
\longrightarrow Ho(*/Es_{\mathbb{V}})_{c}$, pour obtenir 
$$B\mathbb{R}\underline{End}(L\omega) :  
Ho(*/Es_{\mathbb{U}})_{c} \longrightarrow Ho(*/Es_{\mathbb{V}})_{c}.$$ 
 
\begin{thm}\label{t4} 
Pour tout $(X,x) \in Ho(Es_{\mathbb{U}})_{c}$, le $\Delta^{o}$-espace 
$\mathbb{R}\underline{End}(L\omega_{x})$ est un $H_{\infty}$-espace. 
De plus, foncteur $(X,x)  
\mapsto B\mathbb{R}\underline{End}(L\omega_{x})$ 
est isomorphe  \`a l'inclusion naturelle  
$Ho(*/Es_{\mathbb{U}})_{c} \longrightarrow Ho(*/Es_{\mathbb{V}})_{c}$. 
\end{thm} 
 
\textit{Preuve:} De mani\`ere \'equivalente il faut montrer  
que pour tout $(X,x) */Es_{\mathbb{U}}^{cw}$, il existe un 
isomorphisme dans $Ho(Pr\Delta^{o}-SEns_{\mathbb{V}})$, fonctoriel en $(X,x)$ 
$$\mathbb{R}\Omega_{x}X \simeq \mathbb{R}\underline{End}(L\omega_{x}).$$ 
En effet, comme $\mathbb{R}\Omega_{x}X$ est un $H_{\infty}$-espace, il en sera 
de m\^eme de $\mathbb{R}\underline{End}(L\omega_{x})$. \\

Tout d'abord, par le corollaire \ref{c2}, le foncteur  
$(X,x) \mapsto \mathbb{R}\underline{End}(L\omega_{x})$ est isomorphe 
\`a $(X,x) \mapsto \mathbb{R}\underline{End}(L\omega_{S(x)})$,  
o\`u $L\omega_{S(x)} : Fib(S(X)) \longrightarrow Top$ est 
le morphisme fibre au point $S(x) \in S(X)$. Rappelons  
que le foncteur $d : Ho(SGp_{\mathbb{U}}) \longrightarrow  
Ho(*/SEns_{\mathbb{U}})_{c}$ est une \'equivalence de cat\'egories  
(o\`u $Ho(*/SEns_{\mathbb{U}})_{c}$ est la sous-cat\'egorie de $Ho(*/SEns_{\mathbb{U}})$  
form\'ee des objets connexes), dont 
nous avons construit un inverse explicite $G :  
Ho(*/SEns_{\mathbb{U}})_{c} \longrightarrow Ho(SGp_{\mathbb{U}})$ (voir le paragraphe $1.1$).  
Ainsi, pour $(S,s) \in Ho(*/SEns_{\mathbb{U}})$, l'isomorphisme naturel 
$(S,s)\simeq d(G(S,s))$ et le lemme \ref{l1} induisent  
un isomorphisme dans $Ho(*/1-PrCat_{\mathbb{V}})$, fonctoriel en $(S,s)$ 
$$(\mathbb{R}\underline{Hom}(Fib(S),Top),L\omega_{s}) \simeq  
(\mathbb{R}\underline{Hom}(Fib(d(G(S,s))),Top),L\omega_{*}).$$ 
Ainsi, d'apr\`es le th\'eor\`eme \ref{p5},  
on dispose d'un isomorphisme fonctoriel en $(S,s)$ 
$$(\mathbb{R}\underline{Hom}(Fib(S),Top),L\omega_{s})  
\simeq (\mathbb{R}\underline{Hom}(G(S,s)-Top,Top),L\omega_{G(S,s)}),$$ 
o\`u $\omega_{G(S,s)} : G(S,s)-SEns_{\mathbb{U}} \longrightarrow SEns_{\mathbb{U}}$  
est le foncteur d'oubli de la structure de 
$G(S,s)$-module. En composant avec $\mathbb{R}\Omega_{*}$ on trouve un isomorphisme fonctoriel en $(S,s)$ 
$$\mathbb{R}\underline{End}(L\omega_{s})\simeq \mathbb{R}\underline{End}(L\omega_{G(S,s)}).$$ 
En r\'esum\'e, on dispose donc pour tout $(X,x) \in Ho(*/Es_{\mathbb{U}})_{c}$,  
d'un isomorphisme fonctoriel en $(X,x)$ 
$$\mathbb{R}\underline{End}(L\omega_{x})\simeq \mathbb{R}\underline{End}(L\omega_{G(S(X),S(x))}).$$ 
 
D'apr\`es le corollaire \ref{c4}, pour tout  
$G \in SGp_{\mathbb{U}}$, le foncteur $\omega_{G} : G-SEns_{\mathbb{U}} \longrightarrow SEns_{\mathbb{U}}$  
(d'adjoint \`a gauche $F$) v\'erifie les hypoth\`eses du corollaire \ref{c1}. Ceci implique que le foncteur  
$$\begin{array}{ccc} 
Ho(*/SEns_{\mathbb{U}})_{c} & \longrightarrow & Ho(Pr\Delta^{o}-SEns_{\mathbb{V}}) \\ 
(S,s) & \mapsto & \mathbb{R}\underline{End}(L\omega_{G(S,s)}) 
\end{array}$$ 
est en r\'ealit\'e isomorphe au foncteur 
$$\begin{array}{ccc} 
Ho(*/SEns_{\mathbb{U}})_{c} & \longrightarrow & Ho(Pr\Delta^{o}-SEns_{\mathbb{V}}) \\ 
(S,s) & \mapsto & \underline{End}_{G(S,s)}(F(*))^{o}, 
\end{array}$$ 
o\`u $\underline{End}_{G(S,s)}(F(*))^{o}$ est le mono\"{\i}de simplicial oppos\'e \`a celui 
des endomorphismes de l'objet $F(*) \in G(S,s)-SEns_{\mathbb{U}}$. 
Mais par d\'efinition, $F(*)$ est l'ensemble simplicial  
sous-jacent \`a $G(S,s)$, sur lequel $G(S,s)$ op\`ere par 
translations \`a gauche. Or, il est imm\'ediat de v\'erifier  
que lorsque $G \in SGp_{\mathbb{U}}$, le morphisme naturel de mono\"{\i}des simpliciaux  
$$G \longrightarrow \underline{Hom}_{G}(G)^{o},$$ 
donn\'e par l'action de $G$ sur lui-m\^eme \`a droite,  
est un isomorphisme. Ainsi, pour tout $(X,x) \in Ho(*/Es_{\mathbb{U}})_{c}$,  
$\mathbb{R}\underline{End}(L\omega_{G(S(X),S(x))})$ est naturellement  
isomorphe \`a $G(S(X),S(x))$. Ceci implique donc que l'on a  
des isomorphismes de $\Delta^{o}$-espaces, fonctoriels en $(X,x)$ 
$$\mathbb{R}\underline{End}(L\omega_{x})\simeq  
\mathbb{R}\underline{End}(L\omega_{G(S(X),S(x))})\simeq G(S(X),S(x)).$$ 
Comme $G(S(X),S(x))$ est un groupe simplicial, et donc 
un $H_{\infty}$-espace, ceci implique en particulier la premi\`ere assertion du th\'eor\`eme. 
 
D'apr\`es le th\'eor\`eme \ref{t1} et la proposition \ref{p1} le foncteur 
$$\begin{array}{ccc} 
Ho(*/Es_{\mathbb{U}})_{c} & \longrightarrow & Ho(H_{\infty}-SEns_{\mathbb{U}}) \\ 
(X,x) & \mapsto & G(S(X),S(x)) 
\end{array}$$ 
est un inverse de  
$B : Ho(H_{\infty}-SEns_{\mathbb{U}}) \longrightarrow  Ho(*/Es_{\mathbb{U}})_{c}$, et donc 
est naturellemet isomorphe au foncteur  
$\mathbb{R}\Omega_{*} : Ho(*/Es_{\mathbb{U}})_{c} \longrightarrow  
Ho(H_{\infty}-SEns_{\mathbb{U}})$. En d'autres termes,  
il existe des isomorphismes fonctoriels de $H_{\infty}$-espaces $G(S(X),S(x))\simeq \Omega_{*}(X,x)$.  
 
En conlusion, pour tout $(X,x) \in Ho(*/Es_{\mathbb{U}})_{c}$,  
on dispose d'isomorphismes fonctoriels de $H_{\infty}$-espaces  
$$\mathbb{R}\underline{End}(L\omega_{x})\simeq G(S(X),S(x))\simeq \mathbb{R}\Omega_{*}(X,x),$$ 
ce qu'il fallait d\'emontrer. \hfill $\Box$ \\ 
 
\begin{cor}\label{c5} 
Pour tout $(X,x)\in */Es_{\mathbb{U}}^{cw}$, soit $Loc(X)$  
la $1$-cat\'egorie de Segal des champs localement constants sur $X$, et $L\omega_{x}$ son 
foncteur fibre en $x$. Alors 
il existe un isomorphisme naturel dans $Ho(Pr\Delta^{o}-SEns_{\mathbb{V}})$ 
$$\mathbb{R}\Omega_{x}X \simeq \mathbb{R}\underline{End}(L\omega_{x}).$$ 
 
En particulier, il existe des isomorphismes fonctoriels 
$$\pi_{m}(X,x)\simeq \pi_{m-1}(\mathbb{R}\underline{End}(L\omega_{x}),Id)).$$ 
\end{cor} 
 
\textit{Remarques:} 
 
\begin{itemize} 

\item  On peut remarquer que le groupe des automorphismes du foncteur 
d'inclusion $Ho(*/Es_{\mathbb{U}})_{c} \longrightarrow Ho(*/Es_{\mathbb{V}})_{c}$
est trivial. Ainsi, l'isomorphisme du th\'eor\`eme \ref{t4} est en r\'ealit\'e
unique. De m\^eme, il existe un unique isomorphisme fonctoriel
$\Omega_{x}X \simeq \mathbb{R}End(\omega_{x})$.

\item  
Pour tout champ localement constant $F \in Ho(Loc(X))$, on dispose  
d'un morphisme de $H_{\infty}$-espaces \textit{d'\'evaluation en $F$} 
$$\mathbb{R}End(\omega_{x}) \longrightarrow \mathbb{R}End(\omega_{x}(F)).$$ 
A travers l\'equivalence du corollaire \ref{c5}, ce morphism induit un 
morphism bien d\'efini de $H_{\infty}$-espaces 
$$\Omega_{x}X \longrightarrow \mathbb{R}End(\omega_{x}(F)),$$ 
qui peut \^etre compris comme \textit{l'$\infty$-monodromie du champ $F$}.  
Le morphisme $\Omega_{x}X \longrightarrow \mathbb{R}End(\omega_{x})$  
est alors une sorte d'\textit{int\'egration} de tous les morphismes 
$\Omega_{x}X \longrightarrow \mathbb{R}End(\omega_{x}(F))$ lorsque 
$F$ varie dans la $1$-cat\'egorie de Segal $Loc(X)$. C'est pour cela  
qu'il m\'erite peut-\^etre le nom d'\textit{$\infty$-monodromie universelle}. 
 
\item Dans sa lettre \`a L. Breen (dat\'ee du $17/02/1975$, voir \cite{gr}) A. Grothendieck 
\'enonce des r\'esultats proches du corollaire \ref{c5}, et qui utilisent 
une hypoth\'etique notion de \textit{champs en $n$-groupo\"{\i}des localement constants}.
Pour \^etre plus pr\'ecis, citons en quelques lignes (haut de la page $3$). \\

\textit{\dots il faudrait donc expliciter comment un $n$-groupo\"{\i}ode $\mathcal{C}_{n}$ 
se r\'ecup\`ere, \`a $n$-\'equivalence pr\`es, par la connaissance de la $n$-cat\'egorie
$\mathcal{T}_{n}:=(n-\underline{Hom})(\mathcal{C}_{n},((n-1)-Cat))$ des $(n-1)$-syst\`emes
locaux sur $\mathcal{C}_{n}$. On aurait envie de dire que $\mathcal{C}_{n}$ est la cat\'egorie
des "$n$-foncteurs fibres" sur $\mathcal{T}_{n}$, i.e. des $n$-foncteurs
$\mathcal{T}_{n} \longrightarrow ((n-1)-Cat)$ ayant certaines
propri\'et\'es d'exactitude \dots}. \\
  
Nous avons envie d'interpr\'eter ces quelques lignes comme une version $n$-tronqu\'ee et non-point\'ee
de notre formule $X \simeq B\mathbb{R}End(\omega_{x})$. 
 
\end{itemize} 
 
\end{section} 
 
\begin{section}{Conclusions} 
 
La conclusion philosophique que l'on peut tirer du corollaire \ref{c5}  
est que le type d'homotopie d'un $CW$ complexe point\'e $(X,x)$ est 
enti\`erement cod\'e dans la donn\'ee cat\'egorique  
$(\mathbb{R}\underline{Hom}(Loc(X),Top),L\omega_{x})$. Ceci permet alors de plonger la 
th\'eorie de l'homotopie des $CW$ complexes point\'es dans  
celle des $1$-cat\'egories de Segal munies de morphismes vers $Top$. Plus 
pr\'ecis\'ement, on peut montrer que si  
$Top \longrightarrow Top'$ est un mod\`ele fibrant, alors le foncteur naturel 
$$\begin{array}{ccc} 
Ho(*/Es^{cw}_{\mathbb{U}})\simeq Ho(*/Es_{\mathbb{U}}) & \longrightarrow & Ho(1-PrCat/Top') \\ 
(X,x) & \mapsto & (Loc(X),L\omega_{x}) 
\end{array}$$ 
est pleinement fid\`ele lorsqu'on le restreint \`a la  
sous-cat\'egorie de $Ho(*/Es_{\mathbb{U}})$ form\'ee des objets connexes par arcs.  
 
On peut ainsi poser raisonnablement la question de  
l'existence d'une notion de $1$-cat\'egories de Segal galoisiennes (voir 
multi-galoisiennes) g\'en\'eralisant la notion de  
cat\'egories galoisiennes de \cite[$\S V$]{sga1}. Cette notion devrait pouvoir  
\^etre d\'evelopp\'ee \`a l'aide du th\'eor\`eme de  
Beck pour les $1$-cat\'egories de Segal, qui n'est malheureusement 
pour le moment qu'une conjecture (voir \cite[$\S 5.10$]{to}),  
et devrait alors redonner une preuve plus conceptuelle 
du th\'eor\`eme \ref{t4} ne faisant pas intervenir de  
cat\'egories de mod\`eles ferm\'ees mais uniquement des $1$-cat\'egories de Segal.  
Une telle th\'eorie permettrait aussi 
de poss\'eder une nouvelle d\'efinition ainsi qu'une  
nouvelle compr\'ehension du type d'homotopie \'etale d'un sch\'ema (comme 
il est d\'efini dans \cite{am} par exemple),  
ou encore de certaines compl\'etions utilis\'ees en topologie alg\'ebrique  
(compl\'etions profinies, nilpotentes \dots).  
J'esp\`ere revenir sur cette question ult\'erieurement.  
 
Enfin, on peut aussi esp\'erer d\'evelopper un  
formalisme Tannakien pour les $1$-cat\'egories de Segal, qui est  
un analogue lin\'eaire de la th\'eorie Galoisienne, et qui permet de 
remplacer les champs localement constants par d'autres  
coefficients tels des pr\'efaisceaux de complexes localement constants.  
Le lecteur int\'eress\'e par ce type de consid\'erations  
pourra se reporter au manuscript temporaire \cite{to}. 
 
\end{section} 
 
\hspace{10mm} 
 
\textit{Remerciements:} Je remercie L. Breen, A. Hirschowitz, L. Katzarkov, C. Simpson, J. Tapia,  et   
G. Vezzosi pour leurs tr\`es utiles commentaires et remarques au cours de discussions et correspondances sur le sujet.  
 
Je tiens aussi \`a remercier l'institut Max Planck pour son hospitalit\'e et ses conditions de travail exceptionnelles.

\textsf{Bertrand To\"en}: Laboratoire J.A. Dieudonn\'e, Universit\'e de Nice Sophia-Antipolis,  
Parc Valrose, 06108 Nice Cedex 2. e-mail: toen@math.unice.fr
 
\end{document}